\documentclass[VANCOUVER,STIX1COL]{WileyNJD-v2}
\newcommand{\VM}[1]{\mathbf{#1}}

\newcommand{\norm}[1]{\left \| #1 \right \|}
\newcommand{\Ten}[1]{\mathcal{#1}}
\newcommand{\R}{\mathbb{R}}
\newcommand{\C}{\mathbb{C}}
\newlength\figW
\usepackage{tikz}
\usepackage{pgfplots}
  
\usepackage{float}
\usepackage{cool}
\usepackage{subcaption}
\pgfplotsset{compat=1.14}
\usepackage{todonotes}
\usepackage{mathtools}
\usepackage{xcolor}

\makeatletter
\renewcommand*\env@matrix[1][*\c@MaxMatrixCols c]{%
  \hskip -\arraycolsep
  \let\@ifnextchar\new@ifnextchar
  \array{#1}}
\makeatother

\usepackage{moreverb}

\newcommand\BibTeX{{\rmfamily B\kern-.05em \textsc{i\kern-.025em b}\kern-.08em
T\kern-.1667em\lower.7ex\hbox{E}\kern-.125emX}}

\articletype{Article Type}

\received{13 November, 2020}
\revised{<day> <Month>, <year>}
\accepted{<day> <Month>, <year>}

\begin{document}

\title{Subspace method for multiparameter-eigenvalue problems based on tensor-train representations}

\author[1]{Koen Ruymbeek*}

\author[1]{Karl Meerbergen}

\author[1]{Wim Michiels}

\authormark{AUTHOR ONE \textsc{et al}}

\address[1]{\orgdiv{Department of Computer Science}, \orgname{KU Leuven}, \orgaddress{\state{Leuven}, \country{Belgium}}}

%

\corres{*Koen Ruymbeek \email{koen.ruymbeek@kuleuven.be}}

\presentaddress{Koen Ruymbeek, Department of Computer Science, KU Leuven, Celestijnenlaan 200A, 3001 Leuven, Belgium}

\abstract[Abstract]{In this paper we solve $m$-parameter eigenvalue problems ($m$EPs), with $m$ any natural number by representing the problem using Tensor-Trains (TT) and designing a method based on this format. $m$EPs typically arise when separation of variables is applied to separable boundary value problems. Often, methods for solving $m$EP are restricted to $m = 3$, due to the fact that, to the best of our knowledge, no available solvers exist for $m>3$ and reasonable size of the involved matrices. 
In this paper, we prove that computing the eigenvalues of a $m$EP can be recast into computing the eigenvalues of TT-operators. We adapted the algorithm in \cite{Dolgov2014a} for symmetric eigenvalue problems in TT-format to an algorithm for solving generic $m$EPs. This leads to a subspace method whose subspace dimension does not depend on $m$, in contrast to other subspace methods for $m$EPS. This allows us to tackle $m$EPs with $m > 3$ and reasonable size of the matrices. We provide theoretical results and report numerical experiments. The MATLAB code is publicly available.}



\maketitle

\maketitle

\section{Introduction} \label{sect:intro}
In this paper we compute solutions $\left( (\lambda_1, \lambda_2, \hdots, \lambda_m), (\VM x_1 , \VM x_2, \hdots, \VM x_m) \right) $ of the $m$-parameter eigenvalue problem  (mEP)
\begin{equation} \label{eqn:mEP}
\left\{\begin{matrix}
\VM A_{1} \VM x_1 & =  & \lambda_1 \VM B_{11} \VM x_1 + \hdots + \lambda_m \VM B_{1m} \VM x_1 \\
\VM A_{2} \VM x_2 & =  & \lambda_1 \VM B_{21} \VM x_2 + \hdots + \lambda_m \VM B_{2m} \VM x_2 \\
& \vdots \\
\VM A_{m} \VM x_m & = & \lambda_1 \VM B_{m1} \VM x_m + \hdots + \lambda_m \VM B_{mm} \VM x_m 
\end{matrix}\right.
\end{equation}
where $\VM A_i, \VM B_{ij} \in \R^{n_i \times n_i}, i,j = 1, \hdots, m$. We call the $m$-tuple $\left( \lambda_1, \lambda_2, \hdots, \lambda_m \right)$ an eigenvalue-tuple and $\left( \VM x_1, \VM x_2, \hdots, \VM x_m\right)$ the associated eigenvector-tuple. Usually one is interested in a selection of eigenvalues. In this paper, we aim for the eigenvalue-tuples whose $i$th component, $\lambda_i$ has smallest modulus. Without loss of generality we assume that $i = m$. We refer to \cite{Atkinson1972} for a detailed description of multiparameter eigenvalue problems. 

Examples of multiparameter eigenvalue problems can be found in the solution of certain separable boundary value problems, see \cite{Atkinson2011}. Here, a system of $m$ linear ordinary differential equations of the form 
\begin{equation} \label{eqn:diff_equation}
p_{j}\left(x_{j}\right) y_{j}^{\prime \prime}\left(x_{j}\right)+q_{j}\left(x_{j}\right) y_{j}^{\prime}\left(x_{j}\right)+r_{j}\left(x_{j}\right) y_{j}\left(x_{j}\right)=\sum_{\ell=1}^{m} \lambda_{\ell} s_{j \ell}\left(x_{j}\right) y_{j}\left(x_{j}\right), \quad j=1, \ldots, m\end{equation}
is obtained with $\VM x_j \in [a_j, b_j]$ and appropriate boundary conditions. Discretisation of \eqref{eqn:diff_equation} leads to a system of the form \eqref{eqn:mEP}.

There are already available methods for large $n_i, i = 1, 2, \hdots, m$ for the case $m = 2$ \cite{Plestenjak2015}, \cite{Meerbergen2015} and $m=3$ \cite{Hochstenbach2019}, but for $m > 3$, as far as we know, there are no methods yet that work for large $n_i, i = 1, 2, \hdots, m$.  In this paper we present an algorithm based on tensors which fills this gap. 

We can associate with \eqref{eqn:mEP} an equivalent generalised eigenvalue problem. We define the $m \times m$-operator determinants 
\begin{equation} \label{eqn:Delta_0}
\VM \Delta_{0}:=\left|\begin{array}{ccc}
\VM B_{11} & \hdots & \VM B_{1 m} \\
\vdots & & \vdots \\
\VM B_{m 1} & \hdots & \VM B_{m m}
\end{array}\right|_{\otimes}=\sum_{\sigma \in S_{m}} \operatorname{sgn}(\sigma) \VM B_{1 \sigma_{1}} \otimes \VM B_{2 \sigma_{2}} \otimes \hdots \otimes \VM B_{m \sigma_{m}}\end{equation}
with $S_m$ the set of permutations of $\{1,\hdots,m \}$ and $\operatorname{sgn}(\sigma)$ defined as the sign of the permutation $\sigma$. 
Analogously $\VM \Delta_i$ is defined as
\begin{equation} \label{eqn:Delta_i}
\VM \Delta_{i}:=\left|\begin{array}{ccccccc}
\VM B_{11} & \hdots & \VM B_{1, i-1} & \VM A_1 & \VM B_{1, i+1} & \hdots & \VM B_{1 m} \\
\vdots & \vdots & \vdots & \vdots & \vdots & \vdots & \vdots \\
\VM B_{m 1} & \hdots & \VM B_{m, i-1} & \VM A_m & \VM B_{m, i+1} & \hdots & \VM B_{m m} \\
\end{array}\right|_{\otimes}, i = 1, \hdots, m.
\end{equation}
Note that we denote here with $\otimes$ the left-Kronecker product \cite{Regalia1989}.
It can be proven that \eqref{eqn:mEP} is equivalent with the generalised eigenvalue problems
\begin{equation} \label{eqn:Delta}
    \VM \Delta_i \VM x  = \lambda_i \VM \Delta_0 \VM x \quad ,i= 1, \hdots m.
\end{equation}
The dimension of $\VM \Delta_i$ is $n_1 n_2 \hdots n_m$, which soon gets high even for modest $n_j, j = 1,\hdots, m$ and $m$. The eigenvalue $\lambda_i$ is the $i$th element of an eigenvalue -tuple of \eqref{eqn:mEP} and eigenvectors $\VM x$ of \eqref{eqn:Delta} are of the form $\VM x_1 \otimes \VM x_2 \otimes \hdots \otimes \VM x_m$ where the $\VM x_j, j = 1, \hdots, m$ together form an eigenvector-tuple of \eqref{eqn:mEP}. In \cite{Meerbergen2015, Hochstenbach2019} for $m = 2$ resp. $3$, the idea is to repeatedly project each equation on a subspace and to solve the smaller projected problems. The problem for $m > 3$ is that even when all subspace dimensions are as low as $10$, the projected eigenvalue problem already has dimension $10^m$, which may because of its size become a challenging problem. In the setting that we present, let all matrices be of order $n$, the dimension of the involved projected problem is independent of $m$ which offers us a significant advantage from a computational point of view when solving $m$EP for $m > 3$.

Our algorithm uses the same ideas as the algorithm in \cite{Dolgov2014a}. To this end, we prove that the $\VM \Delta$-matrices in \eqref{eqn:Delta_0}, \eqref{eqn:Delta_i} can be written in the Tensor Train-format (TT). 
The aforementioned paper only presented a subspace algorithm for symmetric standard eigenvalue problems (with the matrix written in TT-format) which implies that the smallest eigenvalues can be interpreted via a variational characterisation as the result of an optimisation problem. This approach does not apply to the generalised eigenvalue problem $(\VM \Delta_m, \VM \Delta_0)$. We need to pay extra attention in selecting the wanted eigenvalues of the projected problem. Furthermore we have the additional information that in our problem the eigenvectors are of rank one, which we can exploit to increase efficiency. We show the effectiveness of our algorithm by numerical examples. 


Throughout the paper, we denote vectors by small bold letters and matrices by capital bold letters. With calligraphic letters we denote tensors of order strictly larger than two. Furthermore, we denote by $\VM{I}$ the identity matrix of appropriate dimensions and by the function $\text{diag}(\VM{v})$ a diagonal matrix with the elements of the vector $\VM{v}$ on its diagonal. We call a matrix $\VM X \in \R^{n \times k}$ orthonormal if $\VM X^T \VM X = \VM I_k$.
For slicing a matrix or tensor, we use MATLAB notation. Note that if only one index is left free, we interpret the slice as a column vector. If two indices are left free, the slice is a matrix whose row index corresponds to the first free tensor index and the column to the second free tensor index. For example, let $\Ten G$ be a $n_1 \times n_2 \times n_3$-tensor then $\Ten{G}(i,:,k)$ is considered a column vector and $\Ten{G}(:,j,:)$ is an $n_1 \times n_3$ matrix. If we apply an SVD-decomposition on a matrix, we assume that the singular values are in descending order.

\medskip The structure of this paper is as follows. In Section~\ref{sect:algo_Dolgov} we discuss the used TT-format and we explain the method in \cite{Dolgov2014a} that solves a related problem. We explain and motivate the changes to the method from \cite{Dolgov2014a} in Section \ref{sect:Sel_eigv_conv}. In Section~\ref{sect:NumExp} we show with some numerical experiments the effectiveness of our algorithm and the dependency of our algorithm on $m$ and the size of the matrices in a $m$EP. We close this paper with conclusions. The Matlab code can be downloaded following \href{http://twr.cs.kuleuven.be/research/software/delay-control/mEP/}{this link}.

\section{Calculating eigenvalues of symmetric matrices as TT-operators} \label{sect:algo_Dolgov}
In this paper we propose a subspace method where all matrices and vectors are in Tensor-Train format. To make the paper self contained we first give a short description and a more intuitive characterisation which allows us to prove properties more easily. Then we describe the method proposed in \cite{Dolgov2014a} to calculate the smallest eigenvalues of a symmetric standard eigenvalue problem where the matrix is represented in this format. 

\subsection{The Tensor-Train format}
\begin{definition} \label{def:TT-format}
A tensor $\Ten{Y} \in \C^{n_1 \times \hdots \times n_m}$ of order $m$ is in \emph{TT-format} if  
\begin{equation}
\Ten{Y}(i_1,\hdots,i_m) = \VM G_1(i_1) \hdots \VM G_m(i_m), i_k = 1, \hdots, n_k, k = 1, \hdots m    
\end{equation}
with $\VM G_k(i_k)$ an $r_{k-1} \times r_k$ matrix. The length $r_{k-1}$ and width $r_k$ of these matrices $r_k, k= 0,\hdots, m$ are called the \emph{TT-ranks}. It can immediately be observed that $r_0 = r_m = 1$. We call $r_k, k = 1, \hdots, m$ the \emph{interior} ranks. The matrices $\VM G_k$ can be stored as tensor $\Ten{G}_k \in \C^{r_{k-1} \times n_k \times r_k}, k=1,\hdots,m$ with $\Ten{G}_k(:,i_k,:) = \VM G_k(i_k), i_k = 1, \hdots, n_k$. The tensor $\Ten{G}_k, k = 1, \hdots, m$ is called the TT-core. If it is not clear from the context from which tensor the decomposition is made, we add a superscript indicating the tensor. Each dimension of a tensor is called a mode.
\end{definition}

We associate with a vector $\VM y \in \C^{n_1 n_2 \hdots n_m}$ a tensor $\Ten{Y} \in \C^{n_1 \times n_2 \times \hdots \times n_m}$ such that
$$\operatorname{vec}(\Ten Y) = \VM y.$$
This means that for every $i= 1,2, \hdots (n_1 n_2 \hdots n_m)$, there is a unique $(i_1, i_2, \hdots, i_m)$ such that
$$\VM y(i) = \Ten{Y}( i_1, i_2, \hdots, i_m).$$
We say that a vector $\VM y \in \C^{n_1 n_2 \hdots n_m}$ is in TT-format if its associated tensor $\Ten Y$ is in TT-format. We index the vector $\VM y$ with indices $(i_1,\hdots, i_m)$ of its associated tensor. This means that 
$$\VM y(i) = \Ten Y(i_1, i_2, \hdots, i_m) = \VM G_1(i_1) \hdots \VM G_m(i_m).$$

\begin{definition} \label{def:blockTT}
Matrix $\VM{Y} \in \C^{n_1 n_2 \hdots n_m \times b}$ is in \emph{block-TT format} with index $k$ if
\begin{equation} \label{eqn:blockTT}
    \VM Y(i,i_b) = \VM G_1(i_1) \hdots \VM G_{k-1}(i_{k-1}) \hat{\VM G}_k(i_k, i_b) \VM G_{k+1}(i_{k+1}) \hdots \VM G_m(i_m) 
\end{equation}
for a certain $k = 1, \hdots, m$.
The matrix $ \hat{\VM G}_k(i_k, i_b)$ is stored in a tensor $\Ten{G}_k$ of size $r_{k-1} \times n_k \times b \times r_k$. 
\end{definition}

Definition \ref{def:blockTT} means that we make a TT-decomposition of every column in $\VM Y$ (seen as a tensor) such that all cores except one is in common. At first sight, both representations are not so intuitive. Starting from Definition \ref{def:TT-format}, we can see that
\begin{align*}
    \VM Y(i) & = \VM Y( (i_1, i_2, \hdots, i_m)) \\
    & = \VM G_1(i_1) \VM G_2(i_2) \hdots \VM G_m(i_m) \\
    & = \Ten G_1(i_1, :)^T \Ten G_2(:,i_2,:) \hdots \Ten G_m(:,i_m) \\
    & = \sum_{j_1=1}^{r_1} \sum_{j_2=1}^{r_2} \hdots \sum_{j_{m-1}}^{r_{m-1}} \Ten G_1(i_1, j_1) \Ten G_3(j_1,i_2,j_2) \hdots \Ten G_m(j_{m-1}, i_m).
\end{align*}
From this, it follows that 
$$ \VM Y = \sum_{j_1=1}^{r_1} \sum_{j_2=1}^{r_2} \hdots \sum_{j_{m-1}=1}^{r_{m-1}} \Ten G_1(:, j_1) \otimes \Ten G_2(j_1,:,j_2) \otimes \hdots \otimes \Ten G_m(j_{m-1},:), $$
which makes it more clear how Tensor-Trains work.
For the block TT-format, one can see that 
$$ \VM Y(:,i_b) = \sum_{j_1=1}^{r_1} \sum_{j_2=1}^{r_2} \hdots \sum_{j_{m-1}=1}^{r_{m-1}} \Ten G_1(:, j_1) \otimes \Ten G_2(j_1,:,j_2) \otimes \hdots \otimes \hat{ \Ten G}_k(j_{k-1},:,i_b,j_k) \otimes \hdots \otimes \Ten G_m(j_{m-1},:).$$
In \cite{Dolgov2014a}, the eigenvalues of a symmetric matrix $\VM A \in \R^{n_1 n_2 \hdots n_m \times n_1 n_2 \hdots n_m}$ are calculated using a subspace method. The representation with TT-tensors only makes sense if all operations can be executed within the TT-format. Therefore, also matrix $\VM A$ needs to have a representations such that matrix-vector products with a vector in TT-format and projections on subspaces are implemented efficiently. The solution is to see $\VM A$ as a TT-operator.

Elements of $\VM A$ can be indexed by $2$m-tuples $(i_1, \hdots, i_m, j_1, \hdots, j_m)$ where $(i_1, \hdots, i_m)$ resp. $(j_1, \hdots, j_m)$ are the index for the rows resp. the index for the columns. We say that $\VM A$ is a TT-operator if

\begin{equation} \label{eqn:Delta_ttm}
    \VM A(i_1, \hdots, i_m, j_1, \hdots, j_m) = \VM G^A_1(i_1, j_1) \hdots \VM G^A_m(i_m, j_m)
\end{equation} 
with $\VM G^A_k(i_k, j_k)$ a $r^A_{k-1} \times r^A_k$ matrix. The matrices $\VM G^A_k$ can be stored in a tensor $\Ten{G}^A_k \in \C^{r^A_{k-1} \times n_k \times n_k \times r^A_{k}}$. . For now, we assume that the decomposition is given.
Let $\VM y$ be a vector in TT-format of size $n_1n_2 \hdots n_d$ and consider the product $\VM z = \VM A \VM y$, then 
\begin{align*}
    \VM z(i) & = \Ten Z( i_1, \hdots, i_m) \\
     & = \sum_{j_1, \hdots, j_m} \VM A(i_1, \hdots, i_m, j_1, \hdots, j_m) \Ten Y(j_1, \hdots, j_m) \\
    & = \sum_{j_1, \hdots, j_m} \VM G^A_1(i_1, j_1) \hdots \VM G^A_m(i_m, j_m) \VM G^Y_1(j_1) \hdots \VM G^Y_m(j_m) \\
    & = \sum_{j_1, \hdots, j_m} ( \VM G^A_1(i_1, j_1) \otimes \VM G^Y_1(j_1) ) \hdots ( \VM G^A_1(i_m, j_m) \otimes \VM G^Y_m(j_m) ) \\
    & = \left( \sum_{j_1} \VM G^A_1(i_1, j_1) \otimes \VM G^Y_1(j_1) \right) \left( \sum_{j_2} \VM G^A_2(i_2, j_2) \otimes \VM G^Y_2(j_2) \right) \hdots \left( \sum_{j_m} \VM G^A_m(i_m, j_m) \otimes \VM G^Y_m(j_m) \right).
\end{align*}
By defining $\VM G^Z_k(i_k) =  \sum_{j_k} \VM G^A_k(i_k, j_k) \otimes \VM G^X_k(j_k)$, we can see that the TT-format of the matrix-vector product can directly be derived from the TT-format of its components.

\subsection{Algorithm for symmetric matrices as TT-operators}
The method described here is the method in \cite{Dolgov2014a} for a symmetric matrix $\VM A$ written in TT-format. We discuss it as our method is inspired by \cite{Dolgov2014a}. 
Let $\VM X \in \R^{n_1 n_2 \hdots n_m \times b}$ be the current estimate of the eigenvectors for $b$ eigenvalues in block TT-format with index $k$, so
\begin{equation} \label{eqn:block_TT_X}
    \VM X(:,i_b) = \sum_{j_1=1}^{r_1} \sum_{j_2=1}^{r_2} \hdots \sum_{j_{m-1}=1}^{r_{m-1}} \Ten G_1(:, j_1) \otimes \Ten G_2(j_1,:,j_2) \otimes \hdots \otimes \hat{ \Ten G}_k(j_{k-1},:,i_b,j_k) \otimes \hdots \otimes \Ten G_m(j_{m-1},:).
\end{equation}
The idea behind the method is that we update one mode, leaving the other modes fixed. We start with updating the first mode, then we update the second mode and we continue this process until the $m$th mode. When this is finished, we go back in reversed order towards the first mode.
Using \eqref{eqn:block_TT_X}, this means in practice that for the $k$th mode we search for a solution in the subspace 
\begin{equation} \label{eqn:subsp_Xneqk_full}
    \left \{ \sum_{j_1 = 1}^{r_1}\sum_{j_2 = 1}^{r_2} \hdots \sum_{j_{m-1}=1}^{r_{m-1}} \Ten G_1(:,j_1) \otimes \Ten G_2(j_1,:,j_2) \otimes \hdots \otimes \Ten G_{k-1}(j_{k-2},:,j_{k-1}) \otimes \VM Y  \otimes \Ten G_{k+1}(j_{k},:,j_{k+1}) \otimes \hdots \otimes \Ten G_m(j_{m-1},:) | \VM Y \in \C^{n_k \times r_{k-1} r_k} \right \}.
\end{equation} 
Let $\VM X_{\neq k}$ be an orthonormal basis for \eqref{eqn:subsp_Xneqk_full}, then via the variational characterization of eigenvalues, we solve following optimization problem 
\begin{equation} \label{eqn:var_char}
    \min_{\VM X^X \in \R^{r_{k-1} n_k r_k \times b}}  (\VM X^X_k)^T \VM X_{\neq k}^T \VM A \VM X_{\neq k} \VM X^X_k \text{ subject to }\VM (\VM X^X_k)^T \VM X^X_k =\VM I.
\end{equation}
We discuss now how we can construct a basis for \eqref{eqn:subsp_Xneqk_full}. In order to achieve this, we need the concept of interfaces.
In TT-notation, interfaces are defined as the matrices
\begin{equation} \label{eqn:X_geqk_tt}
\VM X^{>k}(:,(j_{k+1},\hdots, j_m) ) = \VM G_{k+1}(j_{k+1}) \hdots \VM G_m(j_m) 
\end{equation}
\begin{equation} \label{eqn:X_leqk_tt}
\VM X^{<k}((j_1,\hdots, j_{k-1}),:) = \VM G_{1}(j_1) \hdots \VM G_{k-1}(j_{k-1}).
\end{equation}
These matrices are of size $r_k \times (n_{k+1}\hdots n_m)$ resp. $(n_1 \hdots n_{k-1}) \times r_{k-1}$. We observe that $\VM X^{<k}$ is completely defined by the first $k-1$ modes of the decomposition, while $\VM X^{>k}$ is completely defined by the last $m-k$ modes.
When in the next step we need $\VM X^{<k \pm 1}$ and $\VM X^{>k \pm 1}$, we just update the cores.
It can be seen that a basis $\VM X_{\neq k}$ for \eqref{eqn:subsp_Xneqk_full} is 
\begin{equation} \label{eqn:frame}
    \VM X_{\neq k} = \VM X^{<k} \otimes \VM I_{n_k} \otimes (\VM X^{>k})^T.
\end{equation}
In the literature $\VM X_{\neq k}$ is called a \emph{frame matrix}.
We want a condition for the frame matrix to be orthonormal. We make an orthonormal basis by orthonormalizing the interfaces $\VM X^{< k}$ and $\VM X^{> k}$. We introduce \emph{left}-and \emph{right orthonormality}. 

\begin{definition} \cite{Dolgov2014a, Holtz2012}
A TT-core $\Ten G_k$ is \emph{left}- resp. \emph{right orthonormal} if 
$$\sum_{j_k = 1}^{n_k} \VM G_k(j_k)^T \VM G_k(j_k) = \VM I_{r_k} $$
resp.
$$\sum_{j_k = 1}^{n_k} \VM G_k(j_k) \VM G_k(j_k)^T = \VM I_{r_{k-1}}.$$
\end{definition}
 Property \ref{prop:orth} puts a condition on the interfaces such that the frame matrix is orthonormal.
\begin{property} \label{prop:orth}
\cite{Holtz2012, Schollwock2011} If $\VM G_1, \hdots, \VM G_{k-1}$ are left-orthonormal and $\VM G_{k+1}, \hdots, \VM G_{k+1}$ are right-orthonormal, then $\VM X_{\neq k}$ is orthonormal.
\end{property}
In the solvers for \eqref{eqn:var_char}, we only need to know how to compute $\VM X_{\neq k}^T \VM A \VM X_{\neq k}\VM v$ for a vector $\VM v$. In Section 5 of  \cite{Holtz2012} they made an algorithm with complexity $\mathcal{O}( r^3 n^2(r^A)^2)$.

Once we have solved \eqref{eqn:var_char}, the new estimate of the eigenvectors is $\VM X = \VM X_{\neq k} \VM X^X_k$. We want to extract from this a basis $\VM X_{\neq k \pm 1}$ for the subspace used in the next step.
To ease the notation we restrict ourselves to the case $m = 3$. Assume that we have just solved \eqref{eqn:var_char} for $k = 1$ and that we want to form $\VM X_{\neq 2}$.

We first reshape $\VM X^X_1$ to a matrix $\VM X^X_\text{resh}$ of size $n_1 \times r_1 b$, then we perform an SVD-decomposition $\VM X_\text{resh}^X = \VM U \VM S \VM Z^T$ and truncate it to a rank $r$ matrix. This implies that 
\begin{equation} \label{eqn:svd_rank}
    \VM X_1^X(:,i_b) \approx \sum_{i=1}^r \VM U(:,i) \otimes \VM S(i,i) \VM Z( (:,i_b), i) = \sum_{i=1}^r \VM U(:,i) \otimes \hat{ \VM Z}( (:,i_b), i).
\end{equation}
Then we proceed with 
\begin{align}
    \VM X(:,i_b) & = \VM X_{\neq 1} \VM X_1^X(:,i_b) \label{eqn:Xk+1_update} \\ 
    & = (\VM I_{n_1} \otimes (\VM X^{>1})^T ) \VM X_1^X(:,i_b) \nonumber \\
    & \approx \left( \sum_{j_2=1}^{r_2} \VM I_{n_1} \otimes \Ten G_2(:,:,j_2)^T \otimes \Ten G_3(j_2, :) \right) \sum_{i=1}^r \VM U(:,i) \otimes \hat{\VM Z}( (:,i_b), i) \label{eqn:Xk+1_update0} \\
    & = \sum_{i=1}^r \sum_{j_2 = 1}^{r_2} \VM U(:,i) \otimes \sum_{j_1=1}^{r_1} \Ten G_2(j_1,:,j_2) \hat{ \VM Z}( (j_1,i_b), i) \otimes \Ten G_3(j_2,:). \label{eqn:Xk+1_update1}
\end{align}
So by defining $\hat{ \Ten G}_1(:,j_1) = \VM U(:, j_1), \hat{ \Ten G}_2(j_1,:,i_b, j_2) = \sum_{i=1}^{r_1} \Ten G_2(i,:,j_2) \hat{\VM  Z}( (i,i_b),j_1), \hat{ \Ten G}_3(j_2,:) = \VM G_3(j_2,:), $ we have found a block TT-decomposition of $\VM X(:,i_b)$ with index $2$. Furthermore, note that by the orthonormality of the SVD-decomposition, the core $\hat{\VM G}_1$ is left-orthonormal by construction, so $\VM X_{\neq 2}$ is via Property \ref{prop:orth} an orthonormal basis by construction.  This algorithm is stated in Algorithm~ \ref{algo:orig_Dolgov}. In \cite{Dolgov2014a} the author posed the idea to add extra random vectors and to use in \eqref{eqn:Xk+1_update0} and \eqref{eqn:Xk+1_update1}, 
\begin{equation} \label{eqn:kickrank}
\tilde{\VM U} = [\VM U, \VM U_1] \in \R^{n_1 \times (r+p)}, \tilde{\VM Z} = [\hat{ \VM Z}, 0] \in \R^{ r_1 b \times (r+p)} 
\end{equation}
with $\VM U_1$ a random matrix such that the columns of $\tilde{\VM U}$ remain orthonormal against each other, instead of $\VM U$ and $\VM Z$. Adding extra randomness improves the robustness of the algorithm for finding the wanted eigenvalues but makes it also more computationally demanding. 

\begin{algorithm}[h] 
\hspace*{\algorithmicindent} \textbf{Input:} TT-decomposition of $\VM A$ \\
\hspace*{\algorithmicindent} \textbf{Output:}  Eigenpairs of $\VM A$ . 
\begin{algorithmic}[1]
\State Make initial guess $\VM X$ in block-TT format.
\State Orthogonalise $\VM X_{\neq 1}$ by modifying each core using Property \ref{prop:orth}.
\While{ Stopping criterion not fulfilled}
\For{ $k = 1, \hdots, m-1$ (left-to-right half-sweep)}
\State Update $\VM X_{\neq k}$ analogously to \eqref{eqn:Xk+1_update}-\eqref{eqn:Xk+1_update1} (except in the very first iteration) 
\State Solve \eqref{eqn:var_char} and store the solution in $\VM X^X$
\State Reshape $\VM X^X$ to $r_{k-1} \times n_k \times b \times r_k$  
\EndFor
\State Update cores in the right-to-left half-sweep (analogous to left-to-right half-sweep).
\EndWhile
\end{algorithmic}
\caption{Algorithm from \cite{Dolgov2014a} to compute $b$ eigenpairs of a symmetric matrix $\VM A$ in TT-format.}  \label{algo:orig_Dolgov}
\end{algorithm}

We want to use the same ideas to calculate the eigenvalues of the generalized eigenvalue problem in \eqref{eqn:Delta}. We first prove that the $\VM \Delta$-matrices can be interpreted as TT-operator. In contrast to a symmetric eigenvalue problem, the characterization as an optimization problem does no longer hold for the generalized eigenvalue problem
\begin{equation} \label{eqn:proj_problem}
    \left( \VM X_{\neq k}^T \VM \Delta_m \VM X_{\neq k}, \VM X_{\neq k}^T \VM \Delta_0 \VM X_{\neq k} \right).
\end{equation} 
This makes the selection of $b$ eigenpairs of \eqref{eqn:proj_problem} more difficult as it is not clear which eigenvalue to take. As we design the algorithm for solving $m$EPs, we use the induced structure to make it as efficient as possible.

\section{Tensor-algorithm for multiparameter eigenvalue problems} \label{sect:basis}
In this and the following section we describe our algorithm that is implemented for $m$EPs. This section is devoted to the more basic aspects of the algorithm, the next section is more focused on some more technical issues like the convergence criterion and the selection of $b$ eigenvalues of \eqref{eqn:proj_problem}. We start this section by motivating why a method similar to the algorithm in previous section can be applied for solving a $m$EP. The algorithm is given in pseudo-code in Algorithm~\ref{algo:mEP_algo}.

\subsection{$\VM \Delta$-matrices as TT-operators} \label{sect:Delta_rankX}
We have seen in the previous section that the complexity depends on the ranks of the involved matrix $\VM A$ written as TT-operator. In order to use the TT-decomposition in our application, the ranks of the $\VM \Delta$-matrices, expressed as TT-operators must be modest. Without loss of generality we focus on $\VM \Delta_0$. As described in \eqref{eqn:Delta_0}, $\VM \Delta_0$ is a $m \times m$-operator determinant. This means that also each element of $\VM \Delta_0$ is in fact the determinant of a matrix
$$ \VM \Delta_0(i_1, \hdots, i_m, j_1, \hdots, j_m) = \text{det} \left( \begin{bmatrix}
\VM B_{11}(i_1,j_1) & \VM B_{12}(i_1,j_1) & \hdots & \VM B_{1m}(i_1,j_1) \\ 
\VM B_{21}(i_2,j_2) & \VM B_{22}(i_2,j_2) & \hdots & \VM B_{2m}(i_2,j_2) \\ 
\vdots &  \vdots & \vdots & \vdots \\ 
\VM B_{m1}(i_m,j_m) & \VM B_{m2}(i_m,j_m) & \hdots & \VM B_{mm}(i_m,j_m) \\ 
\end{bmatrix} \right), $$
so the cores $\VM G_k^\Delta(i_k,j_k)$ for the operator $\VM \Delta_0$ need to fulfil 
$$ \text{det} \left( \begin{bmatrix}
\VM B_{11}(i_1,j_1) & \VM B_{12}(i_1,j_1) & \hdots & \VM B_{1m}(i_1,j_1) \\ 
\VM B_{21}(i_2,j_2) & \VM B_{22}(i_2,j_2) & \hdots & \VM B_{2m}(i_2,j_2) \\ 
\vdots &  \vdots & \vdots & \vdots \\ 
\VM B_{m1}(i_m,j_m) & \VM B_{m2}(i_m,j_m) & \hdots & \VM B_{mm}(i_m,j_m) \\ 
\end{bmatrix} \right) = \VM G_1^\Delta(i_1,j_1) \hdots \VM G_m^\Delta(i_m,j_m) \quad i_k, j_k = 1, \hdots, n_k, k = 1, \hdots, m.$$
This means we need to find a way to write the determinant as a product of matrices. This can be done using following theorem.

\begin{theorem}  \label{the:determinant} \cite{Jurkat1966} Let $\VM A \in \C^{n \times n}$ and denote with $a_{i:}$ the $i$th row of $\VM A$, then 
\begin{equation}
    \text{det} ( \VM A ) = \VM D_{1,n}(a_{1:}) \VM D_{2,n}(a_{2:}) \hdots \VM D_{n,n}(a_{n:})
\end{equation} 
where $\VM D_{i,n}, i = 1, \hdots, n$ is calculated recursively using following expressions
\begin{align*}
    \VM D_{1,n}( \VM a) & = [\VM a(1), \VM a(2), \hdots, \VM a(n)] \in \C^{1 \times n}\\ 
    \VM D_{n,n}( \VM a) & = [\VM a(n), -\VM a(n-1), \hdots, (-1)^{n-1} \VM a(1)]^T   \in \C^{n  \times 1}\\
    \VM D_{k,n}( \VM a) & = \begin{bmatrix}
    \VM D_{k-1, n-1}( \VM a(2), \hdots, \VM a(n) ) & 0 \\
    (-1)^{k-1} \VM a(1) \VM I & \VM D_{k,n-1} ( \VM a(2), \hdots, \VM a(n))
    \end{bmatrix} \in \C^{  {n \choose k-1} \times {n \choose k} }\quad , k = 2, \hdots, n-1
\end{align*}
\end{theorem}
\begin{corrolary} \label{corr:Delta_rank}
A possible TT-decomposition of $\VM \Delta_0$ is given by 
$$\VM G_k^\Delta(i_k, j_k) = \VM D_{km}( [\VM B_{k1}(i_k,j_k), \hdots, \VM B_{km}(i_k, j_k)]), k= 1, \hdots,m.$$
The ranks $r^\Delta_i, i = 1, \hdots, m$ are described by the $m+1$-st line of Pascal's triangle. Note that these are upper bounds for the internal rank.
\end{corrolary}
Corollary \ref{corr:Delta_rank} gives rise to a recursive algorithm to compute a TT-decomposition of the $\VM \Delta$-matrices. It follows that the ranks are on average $\dfrac{2^m}{m}$ which becomes high if $m$ increases. In Section~\ref{sect:NumExp} we devote a numerical experiment to the influence of $m$ on the computational time. We see that a significant reduction in time can be achieved if we first round off the TT-decomposition upto a certain low tolerance.

The requirement to keep the ranks of the TT-cores small is also the reason why we did not include deflation in Algorithm \ref{algo:mEP_algo}, even though the preservation of the TT-format can be easily achieved by using Hotelling's deflation. We can namely deflate the $q$ already calculated eigenvalues towards infinity by defining 
\begin{equation}
    \VM \Delta_0^{\text{defl}} = \VM \Delta_0 - \sum_{i=1}^q \dfrac{1}{\VM y_i^T \VM \Delta_m \VM x_i} \VM \Delta_0 \VM x_i \VM y_i^T \VM \Delta_m
\end{equation} where $\VM y_i$ is the left-eigenvector (which is also of rank one in TT-format). However it follows that $$\sum_{i=1}^q \dfrac{1}{\VM y_i^T \VM \Delta_m \VM x_i} \VM \Delta_0 \VM x_i \VM y_i^T \VM \Delta_m$$ has internal ranks $r_i + q (r_i^\Delta)^2$ which is too large for practical use. Note that in the context of multiparameter eigenvalue problems, finding a deflation method that respects the TT-rank structure is still an open problem.


\subsection{Size of the projected eigenvalue problem and the rank of $\VM X$} \label{sect:size}
The subspaces $\VM X_{\neq k}$ are of dimension $r_{k-1} n_k r_k$ and so is the projected eigenvalue problem $\left( \VM X_{\neq k}^T \VM \Delta_m \VM X_{\neq k}, \VM X_{\neq k}^T \VM \Delta_0 \VM X_{\neq k} \right)$. This means we can keep the eigenvalue problem small if the internal ranks of $\VM X_{\neq k}$ are low. These ranks are determined by the number of vectors taken in the SVD-decomposition when constructing $\VM X_{\neq k \pm 1}$ from $\VM X_{\neq k}$ and the computed eigenvectors of the projected eigenvalue problem.

When going from $\VM X_{\neq k}$ (see \eqref{eqn:svd_rank}). We know that all eigenvectors of $\left( \VM \Delta_m, \VM \Delta_0 \right)$ are of rank one. For estimating which rank is appropriate, we consider the case where the eigenvectors stored in $\VM X$ have converged. Note that when having a complex-valued eigenvector we split it into its real and complex part as we want to keep the subspace real. A block-TT decomposition of $\VM X$ has internal ranks at most $b$. Let $\VM X$ consist of $p$ real eigenvectors $\VM x^i_1 \otimes \VM x^i_2 \otimes \hdots \otimes \VM x^i_m$ and $\dfrac{p-b}{2}$ complex-valued, then a possible block-TT decomposition of index one of $\VM X$ is 
\begin{align*}
\Ten G_1(:,i_b,i_b) = \VM x_1^{i_b}, \quad \Ten G_2(i_b,:,i_b) = \VM x_2^{i_b}, \quad \VM G_m(i_b,:) = \VM x_m^{i_b}, &  \quad i_b = 1,\hdots,p \\
\Ten G_1(i_b,:,i_b:i_b+1) = [ \text{Re}(\VM x_1^{i_b}),  -\text{Im}(\VM x_1^{i_b}) ], \Ten G_1(:,i_b+1,i_b:i_b+1) = [\text{Im}(\VM x_1^{i_b}), \text{Re}(\VM x_1^{i_b}) ] \\
\Ten G_j(i_b,:,i_b:i_b+1) = [\text{Re}(\VM x_j^{i_b}), -\text{Im}(\VM x_j^{i_b})],
\Ten G_j(i_b+1,:,i_b:i_b+1) = [\text{Im}(\VM x_j^{i_b}), \text{Re}(\VM x_j^{i_b})],& \quad  j = 2, \hdots, m-1 \\
\VM G_m(i_b:i_b+1,:) = [\text{Re}(\VM x_m^{i_b}), \text{Im}(\VM x_m^{i_b})], & \quad i_b = k+1, k+3, \hdots, b-1.
\end{align*}
We consider $r$ to be maximal $b + 1$ as mostly not all eigenvalues are converged at the same time.

\subsection{Stopping criterion} \label{sect:stop_crit}
Usually after a few full sweeps, we find already converged eigenvalues. These are eigenvalues with small modulus, but usually not among the smallest. To make our algorithm more robust in finding the smallest eigenvalues, we do extra full sweeps. We stop after a fixed number of sweeps or when we did not find new eigenvalues in the last sweeps. We noticed that the more clustered the eigenvalues are, the more sweeps it needs to find them all.

\begin{algorithm}[h] 
\hspace*{\algorithmicindent} \textbf{Input:} $m$EP \\
\hspace*{\algorithmicindent} \textbf{Output:} Eigenvalue-tuples $\left( \lambda_1, \lambda_2, \hdots, \lambda_m \right)$ with corresponding eigenvector-tuples $\left( \VM x_1, \hdots, \VM x_m \right)$ which are smallest in absolute value
\begin{algorithmic}[1]
\State Construct $\VM \Delta_m$ and $\VM \Delta_0$
\State Make initial guess $\VM X$ in block-TT format
\State Orthogonalise $\VM X_{\neq 1}$ by modifying each core using Property \ref{prop:orth}
\While{ Stopping criterion not fulfilled}
\For{ $k = 1, \hdots, m-1$ (left-to-right half-sweep)}
\State Update $\VM X_{\neq k}$ using \eqref{eqn:Xk+1_update}-\eqref{eqn:Xk+1_update1} (except in the very first iteration)  \label{line:update}
\State Select $b$ eigenpairs of $\left( \VM X_{\neq k}^T \VM \Delta_m \VM X_{\neq k}, \VM X_{\neq k}^T \VM \Delta_0 \VM X_{\neq k} \right)$ and check their convergence (see Section~\ref{sect:Sel_eigv_conv})  \label{line:sel_conv}
\State Store the associated eigenvectors in $\VM X^X_k$ and make real
\State Reshape $\VM X^X_k$ to $r_{k-1} \times n_k \times b \times r_k$  \label{line:reshape}
\EndFor
\State Update cores in the right-to-left half-sweep (analogous to left-to-right half-sweep)
\EndWhile
\end{algorithmic}
\caption{Algorithm to calculate eigenvalue-tuples with smallest $\lambda_m$ in absolute value of a $m$EP. Line \ref{line:sel_conv} is further discussed in Section~\ref{sect:Sel_eigv_conv}. }  \label{algo:mEP_algo}
\end{algorithm}

\section{Selecting eigenvalues and convergence criterion} \label{sect:Sel_eigv_conv}
In the previous section we showed that it is possible to make an algorithm using the same ideas as in \cite{Dolgov2014a}. There are two main aspects that are completely different: the selection of eigenvalues of the projected problem and the convergence criterion. To achieve this we rely on the fact that eigenvectors of $\left( \VM \Delta_m, \VM \Delta_0 \right)$ are always of rank one in their TT-representation.

\subsection{Selection of eigenvalues} \label{sect:SelectEigenv}
For standard eigenvalue problems with symmetric matrices, calculating the $b$ smallest eigenvalues of $\VM A$ corresponds to minimizing a cost function inferred from the variational characterization of eigenvalues. Projection methods work well \cite{Sirkovic2016a, ruymbeek2019calculating} for symmetric eigenvalue problems as the eigenvalues of the projected eigenvalue problem are bounded by the eigenvalues of the large problem \cite{Parlett1998}. For our non-Hermitian generalised eigenvalue problem 
\begin{equation} \label{eqn:proj_problem1}
\left( \VM X_{\neq k}^T \VM \Delta_m \VM X_{\neq k},\VM X_{\neq k}^T \VM \Delta_0 \VM X_{\neq k} \right)
\end{equation}
 the convergence is much harder to understand. In particular, it is not always clear which Ritz pairs should be selected. Therefore we need to take $b$ significantly high such that if we select an eigenvalue that is not a good estimate of an eigenvalue, it does not harm much the convergence. On the other hand, in Section~\ref{sect:size} we saw that the maximal rank $r$ is chosen as $b+1$ , so as the size of \eqref{eqn:proj_problem1} is $r^2 n$, we cannot take $b$ too large. Via numerical experiments, we concluded that $b$ at least $5$ gives a good balance between the two.

 The following method is heuristic and is based on the method in \cite{Ruymbeek2020}.
\begin{enumerate}
\item Check if the Ritz value is converged (see Section~\ref{sect:Conv_criterion})
\item If it has not converged, we check if the Ritz vector is close enough to one of the Ritz vectors in the previous iteration, see around \eqref{eqn:est_eigenvector}.
\item If there are too many Ritz values that satisfy 2., we select among them the Ritz values with smallest residual norm.
\item If there are not enough Ritz values that fulfil 2., we complete the selection with the Ritz values with smallest residual norms. If there are still not enough, then we add random vectors.
\end{enumerate}

We do not lock Ritz vectors in our subspace neither we deflate them as converged Ritz values as they do not necessarily correspond with the smallest eigenvalues as discussed in Section~\ref{sect:stop_crit}. Let $q$ be the number of already found eigenvalues, then we select the 2b + q Ritz values via the \verb|eigs| (ARPACK) command in Matlab that have smallest modulus and positive real part. Here, one of the two projection matrices need to be formed explicitly which has a cost of $\mathcal{O}(r^4 n^2 (r^\Delta)^2)$, while for the other projection matrix only the matrix-vector product is performed with cost $\mathcal{O}(r^3 n^2 (r^\Delta)^2)$. 

There are two problems with this procedure:
\begin{enumerate}
	\item Within a left-to-right sweep we project the eigenvalue problem on $\VM X_{\neq k}$ in step $k$ and in step $k+1$ we project on $\VM X_{\neq k+1}$. The eigenvectors of the projected eigenvalue problem in step $k$ are not comparable to the eigenvectors of the projected eigenvalue problem in step $k+1$. (an analogous problem arises in right-to-left sweeps)
	\item Calculating the norm of residuals in every step for all $2b + q$ eigenvalues is  expensive.
\end{enumerate}

Without loss of generality we consider the case we update the $k$th mode and we update the $k+1$st mode in the next step. We can solve both problems at the same time. Also for checking the convergence we use the same ideas.
Given a Ritz pair $\left( \mu, \VM X_{\neq k} \VM x^x_k \right)$ of $\left( \VM X_{\neq k}^T \VM \Delta_m \VM X_{\neq k},\VM X_{\neq k}^T \VM \Delta_0 \VM X_{\neq k} \right)$, we make a frame matrix $\hat{\VM X}_{\neq k+1}$ based on only this eigenpair in the same way as in \eqref{eqn:Xk+1_update}-\eqref{eqn:Xk+1_update1}. 
Our estimate $\hat{\VM x}_{k+1}^x$ is then the vectorisation of $\hat{ \Ten G}_{k+1}$. In the next iterate, we compare an eigenvector $\VM x^x_{k+1}$ with $\hat{\VM x}_{k+1}^x$ by  seeing $\VM x^x_{k+1}$ resp. $\hat{\VM x}_{k+1}^x$ as a tensor of size $r_k \times n_{k+1} \times r_{k+1}$ resp. $\hat{r}_k \times n_{k+1} \times r_{k+1}$, make a rank-one TT-approximation 
\begin{equation} \label{eqn:est_eigenvector0}
\VM x^x_{k+1} \approx \VM x^{<k+1} \otimes \overline{\VM x}_{k+1} \otimes \VM x^{> k+1} \nonumber
\end{equation}
resp.
\begin{equation}
\hat{\VM x}_{k+1}^x \approx \hat{ \VM x}^{<k+1} \otimes \hat{\overline{ \VM x}}_{k+1} \otimes \hat{ \VM x}^{> k+1} \label{eqn:est_eigenvector}
\end{equation}
of these tensors and we consider them equal if 
\begin{equation} \label{eqn:similar_vector}
\cos( \theta( \overline{\VM x}_{k+1}, \hat{ \overline{\VM x}}_{k+1})) > 0.99.
\end{equation} 
In words, we check whether the angle between the second mode (i.e. the mode that contains the approximation of one of the vectors in an eigenvector-tuple) of $\hat{\VM x}_{k+1}^x$ and $\VM x_{k+1}^x$ is low enough.
Certainly if convergence starts, from the following Property \ref{prop:xp} we deduce that both vectors can be well approximated by a rank-one tensor.
\begin{property} \label{prop:xp}
Assume that an eigenvector $\VM x = \VM x_1 \otimes \VM x_2 \otimes \hdots \otimes \VM x_m$ of $\left( \VM \Delta_m,  \VM \Delta_0\right)$ is such that it lays in the subspace spanned by the columns of $\VM X_{\neq k}$.
This means that there is a $\VM x^x_k$ such that
$$\VM X_{\neq k} \VM x^x_k = \VM x,$$.
It holds for all $p = 1, \hdots, m$ that there is a frame matrix $\VM X_{\neq p}$ such that
$$\VM x^x_p = \VM X_{\neq p}^T \VM x= \VM x^{<p} \otimes \VM x_p \otimes \VM x^{>p}.$$
The vector $\VM x_p$ is the $p$th element of an eigenvector-tuple of the $m$EP.
\begin{proof}
We have seen that if a vector $\VM x$ lays in the subspace spanned by the columns of $\VM X_{\neq k}$, we can build a frame matrix for mode $k+1$ that also contains this $\VM x$. Therefore, it is enough to prove the result for mode $k$. The latter follows from
\begin{align*}
    \VM x^x_k & = \VM X_{\neq k}^T \left( \VM x_1 \otimes \VM x_2 \otimes \hdots \otimes \VM x_m \right) \\
    & = \left( (\VM X^{<k})^T \otimes \VM I_k \otimes (\VM X^{>k}) \right)( \VM x_1 \otimes \hdots \otimes \VM x_m) \\
    & = \sum_{j_1 = 1}^{r_1} \sum_{j_2 = 1}^{r_2} \hdots \sum_{j_{k-2} = 1}^{r_{k-2}}  \sum_{j_{k+2} = 1}^{r_{k+2}} \hdots \sum_{j_{m-1} = 1}^{r_{m-1}}  \VM G_1(:,j_1)^T \VM x_1 \otimes  \Ten G_2(j_1,:,j_2)^T \VM x_2 \hdots \otimes  \Ten G_{k-2}(j_{k-3}:,j_{k-2})^T \VM x_{k-2}  \otimes \Ten G_{k-1}(j_{k-2},:,:)^T \VM x_{k-1} \\
    & \quad \otimes \VM x_k \otimes \left( \Ten G_{k+1}(:,:,j_{k+1})^T \VM x_{k+1} \otimes \Ten G_{k+2}(j_{k+1},:,j_{k+2})^T \VM x_{k+2} \hdots \otimes \Ten G_{m-1}(j_{m-2},:,j_{m-1})^T \VM x_{m-1} \otimes \VM G_m(j_{m-1}, :)^T \VM x_m \right) \\
    & = \underbrace{ \sum_{j_{k-2}=1}^{r_{k-2}} a_{j_{k-2}}^{k-2} \Ten G_{k-1}(j_{k-2},:,:)^T \VM x_{k-1}}_{ := \VM x^{<k}} \otimes  \VM x_k \otimes \underbrace{\sum_{j_{k+1} = 1}^{r_{k+1}}  a_{j_{k+1}}^{k+1} \Ten G_{k+1}(:,:,j_{k+1}) \VM x_{k+1} }_{:=\VM x^{>k}}
\end{align*}
where the constants $a^{k-2}_{j_{k-2}}$ and $a^{k+1}_{j_{k+1}}$ are obtained from the projections of $\Ten G_{j_i}(j_{i-1},:,j_i)$ on $\VM x_i$ for $i = 1, \hdots, k-2$ resp. $i = k+1, \hdots, m$.
\end{proof}
\end{property}
To make an estimation of the residual norm for eigenpair $(\mu, \VM x^x_k)$ of the projected problem, we use the same frame matrix $\hat{ \VM X}_{\neq k+1}$ and we compute the norm of the residual projected on the subspace spanned by the columns of $\hat{ \VM X}_{\neq k+1}$, i.e
\begin{equation} \label{eqn:residual_est}
     \norm{ \hat{ \VM r}} = \norm{ \hat{\VM X}_{\neq k + 1}^T \VM \Delta_m \hat{\VM X}_{\neq k + 1} \hat{ \VM x}^x_{k + 1} - \mu_i \hat{\VM X}_{\neq k + 1}^T \VM \Delta_0 \hat{\VM X}_{\neq k + 1} \hat{ \VM x}^x_{k + 1}},
\end{equation}
where $\hat{\VM X}_{\neq k + 1}^T \VM \Delta_0 \hat{\VM X}_{\neq k + 1} \hat{ \VM x}^x_{k + 1}$ and $\hat{\VM X}_{\neq k + 1}^T \VM \Delta_m \hat{\VM X}_{\neq k + 1} \hat{ \VM x}^x_{k + 1}$ are directly computed using the TT-representations. This is a much cheaper alternative than explicitly calculating the residual norm.  As \eqref{eqn:residual_est} is always smaller than the residual norm of the full size problem, there is no guarantee that if \eqref{eqn:residual_est} is small it also holds for the residual norm. But if this is already high, we know that the residual norm is even higher and it is for sure not a converged eigenvalue. The calculation of \eqref{eqn:residual_est} is used to decide which eigenpairs to select if the criterion in \eqref{eqn:similar_vector} is not restrictive enough and to decide whether it is a possible converged eigenpair, see Section \ref{sect:Conv_criterion}. As we need to do this in every step of a sweep for all $2b+q$ eigenpairs, it would be too expensive to calculate the residual norm for each eigenpair.

The algorithm is given in pseudo-code in Algorithm ~\ref{algo:sel_eigv}.


\begin{algorithm}[h]
\hspace*{\algorithmicindent} \textbf{Input:} 
\begin{enumerate}
\item Frame matrix $\VM X_{\neq k}, k = 1, \hdots, m$ and direction $d = 1$ or $-1$ (left-to-right half sweep or right-to-left half sweep)
\item $\VM \Delta$-matrices $\VM \Delta_m, \VM \Delta_0$
\item $2b + q$ eigenpairs of \eqref{eqn:proj_problem}, with $q$ the number of already found eigenvalues
\item The list $\hat{\VM X}_\text{selected}^k$ with estimates in \eqref{eqn:est_eigenvector}, selected in the previous iterate. 
\end{enumerate} 
\hspace*{\algorithmicindent} \textbf{Output:} 
\begin{enumerate}
\item Selection of $b$ eigenpairs of \eqref{eqn:proj_problem}. This means we store its eigenvectors in the columns of the matrix $\VM X^X_k$, make a rank-one decomposition of $\hat{ \VM x}_{k+ d}^x$ and add $\hat{ \overline{ \VM x}}_{k+ d}$ to $\hat{\VM X}^{k+ d}_\text{select}$ (see \eqref{eqn:est_eigenvector}) .
\end{enumerate}
\begin{algorithmic}[1]
\For{$i=1$ to $2b + q$}
\State Select the $i$th eigenpair $(\mu,  \VM x^x_k)$ of \eqref{eqn:proj_problem}
\State Make subspace $\hat{ \VM X}_{\neq k + d}$ and vector $\hat{\VM x}^x_{k + d}$ (See 
close to \eqref{eqn:est_eigenvector0})
\State Estimate the residual norm \eqref{eqn:residual_est} and check the convergence of $\VM x^x_k$  (see Section~\ref{sect:Conv_criterion}) \label{line:check_conv}
\State Make a rank-one decomposition \eqref{eqn:est_eigenvector} of $\VM x^x_k$ and retrieve $\overline{ \VM x}_k \in \C^{n_k}$
\If{ $\operatorname{max}_{ \VM x \in \hat{\VM X}_{\text{select}}^k}( \cos( \theta(\VM x, \overline{\VM x}_k)) > 0.99$ (check if vectors are close to each other) and not yet converged} 
\State Select this eigenpair 
\EndIf
\EndFor 
\If{we select not enough eigenvectors}
\State Select eigenvectors that were not selected neither converged with the smallest residual
\EndIf
\If{we select too much eigenvectors}
\State Select among the vectors that we selected the eigenvectors with the smallest residual
\EndIf
\end{algorithmic}
\caption{Algorithm to select $b$ eigenvalues among the $2b + q$ eigenvalues with positive imaginary part that we have calculated from \eqref{eqn:proj_problem1}. This algorithm is used on line \ref{line:sel_conv}  of Algorithm~\ref{algo:mEP_algo}. The method for line \ref{line:check_conv} is explained in Section~\ref{sect:Conv_criterion}.}  \label{algo:sel_eigv}
\end{algorithm}

\subsection{Convergence criterion} \label{sect:Conv_criterion}
We use the same convergence criterion as used in \cite{Hochstenbach2019} for $m = 3$,
\begin{equation} \label{eqn:convergence}
\norm{ \begin{bmatrix}
\VM ( \VM A_1 - \lambda_1 \VM B_{11} - \hdots - \lambda_m \VM B_{1m} ) \VM x_1 \\ 
\VM ( \VM A_2 - \lambda_1 \VM B_{21} - \hdots - \lambda_m \VM B_{2m} ) \VM x_2 \\ 
\vdots \\ 
\VM ( \VM A_m - \lambda_1 \VM B_{m1} - \hdots - \lambda_m \VM B_{mm} ) \VM x_m \\
\end{bmatrix} }_\infty< \epsilon
\end{equation}
so an eigenvalue is considered converged if the residual norm for all equations together is below a tolerance $\epsilon$. We choose $\epsilon$ equal to $10^{-6}$ as this is the tolerance that was used in the aforementioned paper. 

Without loss of generality, we only consider the case where we just have updated the $k$th mode and will update the $k+1$st mode in the next step. In Section~\ref{sect:SelectEigenv}, we explained how we estimate the residual norm for $\VM x^x_k$ by making a frame matrix $\hat{\VM X}_{\neq k + 1}$ and to compute the residual projected on the subspace spanned by the columns of this frame matrix (see \eqref{eqn:residual_est}). The idea is now that if this residual norm is below $\epsilon_1$, to continue this process and to construct the frame matrix $\hat{\VM X}_{\neq k + 2}$, estimate $\hat{\VM x}_{k + 2} \in \C^{n_{k + 2}}$ and check whether the residual projected on the subspace spanned by the columns of this frame matrix is below $\epsilon_1$ too. We can continue this process until all modes are visited or till the moment we find a mode for which this approximation of the residual norm is too high. If we reach $m$ and not all modes are visited, then we start again with $k$ and go back in reverse order to mode one. As we always check the residual whithin a subspace, $\epsilon_1$ is chosen low. We use in our algorithm $\epsilon_1 = 10^{-8}$.  If the residual projected on all subspaces $\hat{\VM X}_{\neq 1}, \hat{\VM X}_{\neq 2}, \hdots, \hat{\VM X}_{\neq m}$ is below $\epsilon_1$, we may have found a  converged eigenvalue. In this case we have an estimate for the vectors $\hat{\VM x}_1, \hdots, \hat{\VM x}_m$ with eigenvalue $\hat{\lambda}_m$. In order to check the criterion \eqref{eqn:convergence}, we need associated values $\hat{\lambda}_1, \hdots, \hat{\lambda}_{m-1}$ to complete the eigenvalue-tuple. In \cite{Hochstenbach2019, Plestenjak2000} a Tensor Rayleigh Quotient iteration (TRQI) was described to refine $\hat{\lambda}_m$ and to give associated values $\hat{\lambda}_1, \hdots, \hat{\lambda}_{m-1}.$

If we have found a new eigenvalue, it is still possible that it is an eigenvalue we already found. In \cite{Hochstenbach2019} a method was described to check this for 3EP. We use an analogous tool for the general $m$EP case. A new eigenvector $\hat{\VM x} = \hat{\VM x}_1 \otimes \hat{\VM x}_2 \otimes  \hdots \otimes \hat{\VM x}_m$ is accepted if 
\begin{equation} \label{eqn:new_eig}
    \max_{p=1\hdots,q} \dfrac{(\VM y^p)^* \VM \Delta_0 \hat{\VM x}}{(\VM y^p)^* \VM \Delta_0 \VM x^p} < \xi
\end{equation}
with $\VM x^p = \VM x_1^p \otimes \VM x_2^p \otimes \hdots \otimes \VM x_m^p$ resp. 
$\VM y^p = \VM y_1^p \otimes \VM y_2^p \otimes \hdots \otimes \VM y_m^p$ already found right- and left eigenvectors. For calculating the right eigenvector-tuple $\left( \VM y_1^p, \hdots, \VM y_m^p \right)$ we use the same tool as used in \cite{Hochstenbach2019}. The tolerance $\xi$ is usually chosen as $10^{-4}$. Note that the expression in \eqref{eqn:new_eig} can be calculated by considering $\VM x, \VM y$ and $\VM \Delta_0$ in TT-format. Remark that in theory $(\VM y^p)^* \VM \Delta_0 \hat{ \VM x} = 0$ if $\hat{ \VM x}$ and $\VM y^p$ are associated to a different eigenvalue. Note that this is a different criterion than in \eqref{eqn:similar_vector}. In \eqref{eqn:similar_vector} we only have estimations of a solution of the $k+1$st equation in the $m$EP \eqref{eqn:mEP}. As this is not the solution of a generalized or standard eigenvalue problem, it is not possible to use a similar criterion as in \eqref{eqn:new_eig}.

The implemented algorithm based on the ingredients described above, can be found in Algorithm~\ref{algo:convergence}. This algorithm is the algorithm to which we refer in Algorithm~\ref{algo:sel_eigv} line~\ref{line:check_conv}.

\begin{algorithm}[h] 
\hspace*{\algorithmicindent} \textbf{Input:} 
\begin{enumerate}
\item Direction $d = 1$ (left-to-right sweep), $d = -1$ (right-to-left sweep) 
\item Subspace $\hat{ \VM X}_{\neq k+d}$ 
\item Vector $\hat{ \VM x}^x_{k+d}$ (see line \ref{line:check_conv} in Algorithm~ \ref{algo:sel_eigv})
\end{enumerate}
\hspace*{\algorithmicindent} \textbf{Output:} Estimate of residual $ \hat{\VM r}$ and check if eigenvalue has converged.
\begin{algorithmic}[1]
\State Compute $\norm{ \hat{\VM r}}$ as in \eqref{eqn:residual_est}
\State Decompose $\hat{\VM x}^x_{k+d}$ and compute estimate $\hat{\VM x}_{k+d} \in \C^{n_{k_1}}$
\State $k_1 = k + d$;
\While{ $\norm{ \hat{\VM r}} < \epsilon_1$ and we do not have estimates for all modes}
    \If{$k_1 == m$ or $k_1 == 1$} (If we reach the first/last mode, we need to change direction)
        \State $k_1 = k$; $d = -d$
    \EndIf
    \State $k_1 = k_1 + d$
    \State Make $\hat{\VM X}_{\neq k_1}$ and vector $\hat{\VM x}_{k_1}$
	\State Decompose $\hat{ \VM x}_{k_1}$ and compute estimate $\hat{\VM x}_{k_1} \in \C^{n_{k_1}}$ 
	\State Compute $\norm{ \hat{ \VM r}}$ as in \eqref{eqn:residual_est}
\EndWhile
\If{ All $\norm{ \hat{\VM r}} < \epsilon_1$}
    \State Check if $\hat{\VM x}_1, \hdots, \hat{\VM x}_m$ is an eigenvector-tuple that is not found yet via \eqref{eqn:new_eig} \label{line:nbr_conv_eigenv}
    \If{Condition \eqref{eqn:convergence} is satisfied and the Ritz value is one of the $p$ smallest eigenvalues we found} 
        \State Add the Ritz pair to the list of found eigenpairs
    \EndIf
\EndIf
\end{algorithmic}
\caption{Algorithm to check if a Ritz value has converged. This algorithm is the algorithm to which we refer in Algorithm~\ref{algo:sel_eigv} line~\ref{line:check_conv}.  In the experiments we choose to keep the $p$ found eigenvalues that best fulfill our criterion. We choose in the numerical experiments $p = 4 b$.}  \label{algo:convergence}
\end{algorithm}

\section{Numerical experiments} \label{sect:NumExp}
In this section we test our algorithm on some of the examples discussed in \cite{Hochstenbach2019} to prove the algorithm works. The second goal of this section is to see how the computational time depend on $m$ and $n$ (the size of the involved matrices). The code for Algorithm~\ref{algo:orig_Dolgov} which was the basis for our algorithm, is the code \verb|dmrg_eig| available in the TT-Toolbox \cite{Dolgov2020}. All algorithms written in the previous sections are implemented in Matlab version R2020a. All experiments are performed on an Intel i5-6300U with 2.5 GHz and 8 GB RAM. We choose $b = 5$ and we add in every step of a sweep one extra random vector as explained in \eqref{eqn:kickrank}.

\subsection{Performance}
\begin{figure}[h] 
\setlength{\figW}{6cm} 
 \begin{subfigure}[b]{0.5\textwidth}
 	\hspace{0.7cm}
%
%
%
\definecolor{mycolor1}{rgb}{0.00000,0.44700,0.74100}%
\begin{tikzpicture}

\begin{axis}[%
width=0.951\figW,
height=0.75\figW,
at={(0\figW,0\figW)},
scale only axis,
xmin=    0,
xmax=   20,
ymin=    0,
ymax=   35,
ylabel style={font=\color{white!15!black}},
ylabel={$\lambda_m$},
axis background/.style={fill=white},
legend style={at={(0.03,0.97)}, anchor=north west, legend cell align=left, align=left, draw=white!15!black, only marks}
]
\addplot [color=mycolor1, draw=none, mark=asterisk, mark options={solid, mycolor1}]
  table[row sep=crcr]{%
    1	2.40498\\
    2	5.59867\\
    3	7.46473\\
    4	10.4954\\
    5	12.3327\\
    6	12.6541\\
    7	15.9169\\
    8	16.9754\\
    9	18.4634\\
   10	18.9359\\
   11	21.2369\\
   12	22.594\\
   13	24.7989\\
   14	25.8319\\
   15	27.4539\\
   16	27.5414\\
   17	28.547\\
   18	30.8625\\
   19	31.6355\\
   20	33.0522\\
};
\addlegendentry{$\lambda_m$}

\addplot [color=red, draw=none, mark=o, mark options={solid, red}]
  table[row sep=crcr]{%
    1	2.40498\\
    2	5.59867\\
    3	7.46473\\
    4	10.4954\\
    5	12.3327\\
    6	12.6541\\
    7	15.9169\\
    8	16.9754\\
    9	18.4634\\
   10	18.9359\\
   11	21.2369\\
   12	22.594\\
   13	24.7989\\
   14	25.8319\\
   17	28.547\\
};
\addlegendentry{found $\lambda_m$}

\end{axis}
\end{tikzpicture}%
	\caption{}
\end{subfigure}
\begin{subfigure}[b]{0.5\textwidth}
	\hspace{0.7cm}  
%
%
%
\definecolor{mycolor1}{rgb}{0.00000,0.44700,0.74100}%
\begin{tikzpicture}

\begin{axis}[%
width=0.951\figW,
height=0.75\figW,
at={(0\figW,0\figW)},
scale only axis,
xmin=    0,
xmax=   20,
ymin=  195,
ymax=  205,
ylabel style={font=\color{white!15!black}},
ylabel={$\lambda_m$},
axis background/.style={fill=white},
legend style={at={(0.03,0.97)}, anchor=north west, legend cell align=left, align=left, draw=white!15!black, only marks}
]
\addplot [color=mycolor1, draw=none, mark=asterisk, mark options={solid, mycolor1}]
  table[row sep=crcr]{%
    1	200.606\\
    2	199.275\\
    3	201.032\\
    4	198.449\\
    5	197.997\\
    6	202.052\\
    7	197.919\\
    8	197.668\\
    9	202.466\\
   10	196.558\\
   11	203.466\\
   12	203.498\\
   13	203.546\\
   14	204.004\\
   15	204.07\\
   16	195.916\\
   17	204.174\\
   18	204.578\\
   19	204.593\\
   20	195.324\\
};
\addlegendentry{$\lambda_m$}

\addplot [color=red, draw=none, mark=o, mark options={solid, red}]
  table[row sep=crcr]{%
    2	199.275\\
    4	198.449\\
    6	202.052\\
    7	197.919\\
    8	197.668\\
   11	203.466\\
   12	203.498\\
   14	204.004\\
   16	195.916\\
};
\addlegendentry{found $\lambda_m$}

\end{axis}
\end{tikzpicture}%
	\caption{}
\end{subfigure}
\caption{We plot the wanted eigenvalues together with the eigenvalues that we found for Example \ref{ex:ellipsoid}. In (a) the eigenvalues are well separated and we we can find the first $2b$ eigenvalues. In (b) the eigenvalues are interior and we see the algorithm has much more difficulties in finding the wanted eigenvalues. We note that different runs can give diverent converged eigenvalues as the algorithm is subject to random numbers. }\label{fig:ellipsoidwave}
\end{figure}
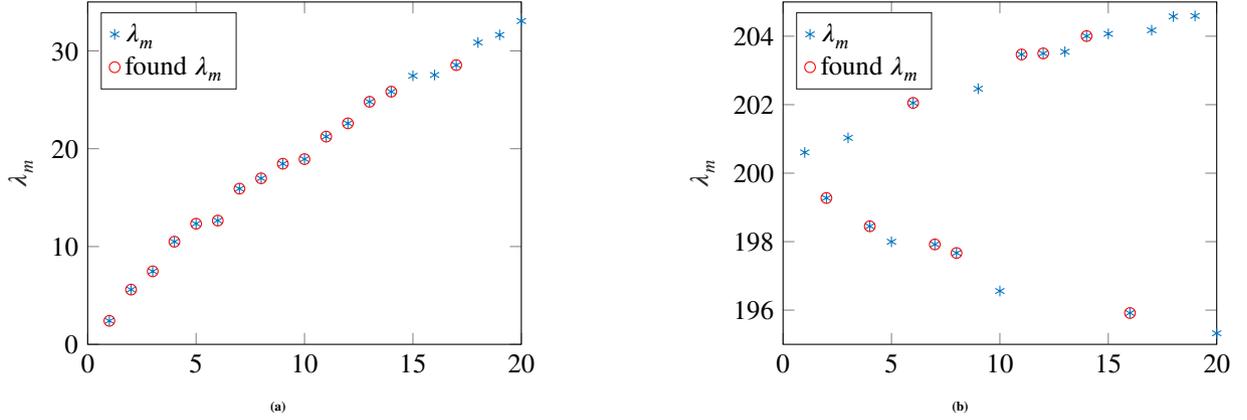

\begin{eexample} \label{ex:ellipsoid}
We use the example described in \cite[Section 5.1]{Hochstenbach2019} for $m = 3$. The only difference is that we use $200$ Chebyshev collocation points instead of 300. We search for the eigenvalues such that $\lambda_3$ is closest to $0$ and $\lambda_3$ closest to $200.$ If the searched eigenvalues are exterior, usually, it is sufficient to do only $20$ full sweeps. For computing interior eigenvalues, we do $100$ full sweeps. In Figure \ref{fig:ellipsoidwave} we show the results. In the first case we found almost the $3b$ eigenvalue-tuples with smallest $\lambda_3$. In the second case, we found almost half of the first $20$ eigenvalues. We note that the algorithm in \cite[Section 5.1]{Hochstenbach2019} can find more eigenvalues for an equal amount of time, but the method \cite{Hochstenbach2019} was especially designed for 3EP and the extension to $m > 3$ is not feasible as was discussed in the introduction (Section~\ref{sect:intro}). 
\end{eexample}

In the literature there are, to the best of our knowledge, no examples for $m > 3$ and reasonable size of the matrices. Therefore, we generate some random examples to test our algorithm. To be able to test how many of the searched eigenvalues we have found, we need to be able to compute all eigenvalues. The technique described here is also the way random examples were created in \cite{Hochstenbach2019}.

We make first real vectors $\VM a_i, \VM b_{ij}, i,j = 1, \hdots, m$ and random linear independent $n \times n$-matrices $\VM U_i, \VM Z_i, i = 1, \hdots, m$ and construct
$$ \VM A_i = \VM U_i \text{diag}( \VM a_i) \VM Z_i$$
$$ \VM B_{ij} = \VM U_i \text{diag}( \VM b_{ij}) \VM Z_i.$$
An eigenvalue-tuple $\lambda = \left( \lambda_1, \hdots, \lambda_m \right)$  of the associated $m$EP is then found as the solution of the linear systems
\begin{equation} \label{eqn:system_random}
\begin{bmatrix}
  \VM b_{11}(i_1) & \VM b_{12}(i_1) & \hdots & \VM b_{1m}(i_1)\\ 
\VM b_{21}(i_2) & \VM b_{22}(i_2) & \hdots & \VM b_{2m}(i_2)\\
\vdots &  & \ddots  & \vdots \\ 
\VM b_{m1}(i_m) & \VM b_{m2}(i_m) & \hdots & \VM b_{mm}(i_m)
\end{bmatrix} \begin{bmatrix} \lambda_1 \\ \lambda_2 \\ \vdots \\ \lambda_m \end{bmatrix} = \begin{bmatrix}
\VM a_1(i_1) \\
\VM a_2(i_2) \\
\vdots \\
\VM a_m(i_m) 
\end{bmatrix}, i_1, \hdots, i_m = 1, \hdots, n.  
\end{equation} 
This means that if we we want to be sure that we compute the $b$ eigenvalue-tuples with smallest $| \lambda_m|$ we need to solve $n^m$ linear systems.

Here we choose to randomly generate $\VM b_{i2}$ and to define $\VM b_{ij} =\VM b_{i2}^{j-1}$. The following code is used to construct $\VM a_i, \VM b_{i2}, \VM U_i$ and $\VM Z_i, i = 1, \hdots, m$. 

\begin{verbatim}
U_c = cell(m,1); V_c = cell(m,1);
B_c = cell(m, m); A_c = cell(m, 1);
for i = 1:m
    U_c{i} = 0.3*rand(n,n)+eye(n);
    V_c{i} = 0.3*rand(n,n)+eye(n);
end
[x, ~] = chebdif_mp( n, 2, 'double');
limit = linspace( -1.9, 2, 2*m+1);
limit = limit(1:2*m);
a_m = -5*randn(n,m); b_m = zeros(n,m);
for i = 1:m
    b_m(:,i) = (x/2*(limit(2*i)-limit(2*i-1))+1/2*(limit(2*i-1)+limit(2*i)));
end
for i = 1:m
    A_c{i} = V_c{i}*diag(a_m(:,i))*U_c{i};
    for j = 1:m
        B_c{i,j+1} = V_c{i}*diag( b_m(:,i).^(j-1))*U_c{i};
    end
end
\end{verbatim}
where \verb|chebdiff_mp| can be found in the package \verb|MultiParEig| \cite{Plestenjak2020}. This code is inspired on the construction of the matrices in Example \ref{ex:ellipsoid}. We can shift the $m$EP such that all $\lambda_m$ are positive-valued by defining a new $m$EP with matrices 

\begin{align}
\hat{ \VM A}_i & := \VM A_i + \eta \VM B_{im}, \quad \hat{ \VM B}_{ij} :=  \VM B_{ij}, \quad i,j = 1, \hdots,m  \label{eqn:shift_eqns}
\end{align}
for a certain shift $\eta \in \R$. It can be seen that this means in the new $m$EP only $\lambda_m$ is shifted to $\lambda_m + \eta$.

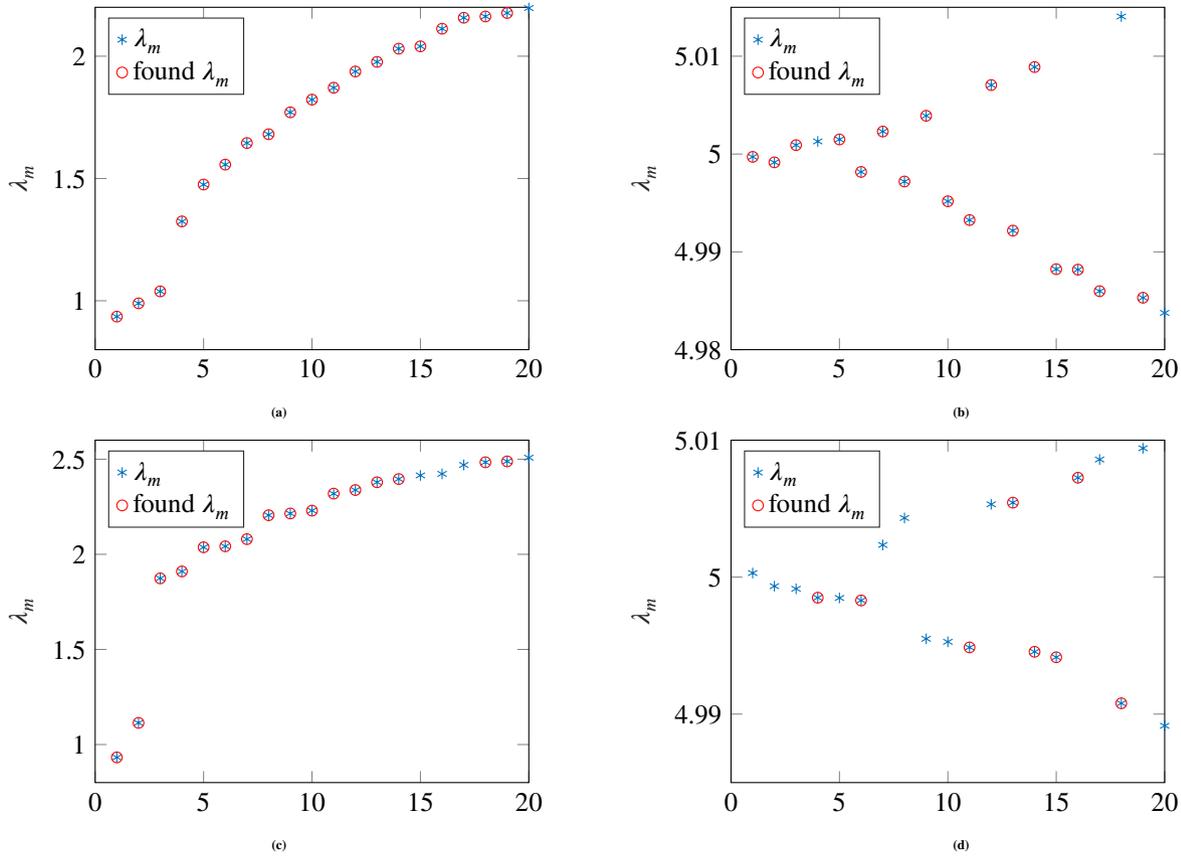
\begin{figure}[h] 
\setlength{\figW}{6cm} 
 \begin{subfigure}[b]{0.5\textwidth}
 	\hspace{0.7cm}
%
%
%
\definecolor{mycolor1}{rgb}{0.00000,0.44700,0.74100}%
\begin{tikzpicture}

\begin{axis}[%
width=0.951\figW,
height=0.75\figW,
at={(0\figW,0\figW)},
scale only axis,
xmin=    0,
xmax=   20,
ymin=  0.8,
ymax=  2.2,
ylabel style={font=\color{white!15!black}},
ylabel={$\lambda_m$},
axis background/.style={fill=white},
legend style={at={(0.03,0.97)}, anchor=north west, legend cell align=left, align=left, draw=white!15!black, only marks}
]
\addplot [color=mycolor1, draw=none, mark=asterisk, mark options={solid, mycolor1}]
  table[row sep=crcr]{%
    1	0.935212\\
    2	0.989598\\
    3	1.03831\\
    4	1.32443\\
    5	1.47521\\
    6	1.55673\\
    7	1.64494\\
    8	1.68094\\
    9	1.77055\\
   10	1.82234\\
   11	1.87107\\
   12	1.9374\\
   13	1.97672\\
   14	2.03069\\
   15	2.04024\\
   16	2.11233\\
   17	2.15719\\
   18	2.16236\\
   19	2.17657\\
   20	2.19681\\
};
\addlegendentry{$\lambda_m$}

\addplot [color=red, draw=none, mark=o, mark options={solid, red}]
  table[row sep=crcr]{%
    1	0.935212\\
    2	0.989598\\
    3	1.03831\\
    4	1.32443\\
    5	1.47521\\
    6	1.55673\\
    7	1.64494\\
    8	1.68094\\
    9	1.77055\\
   10	1.82234\\
   11	1.87107\\
   12	1.9374\\
   13	1.97672\\
   14	2.03069\\
   15	2.04024\\
   16	2.11233\\
   17	2.15719\\
   18	2.16236\\
   19	2.17657\\
};
\addlegendentry{found $\lambda_m$}

\end{axis}
\end{tikzpicture}%
	\caption{}
\end{subfigure}
\begin{subfigure}[b]{0.5\textwidth}
	\hspace{0.7cm}  \\
%
%
%
\definecolor{mycolor1}{rgb}{0.00000,0.44700,0.74100}%
\begin{tikzpicture}

\begin{axis}[%
width=0.951\figW,
height=0.75\figW,
at={(0\figW,0\figW)},
scale only axis,
xmin=    0,
xmax=   20,
ymin= 4.98,
ymax=5.015,
ylabel style={font=\color{white!15!black}},
ylabel={$\lambda_m$},
axis background/.style={fill=white},
legend style={at={(0.03,0.97)}, anchor=north west, legend cell align=left, align=left, draw=white!15!black, only marks}
]
\addplot [color=mycolor1, draw=none, mark=asterisk, mark options={solid, mycolor1}]
  table[row sep=crcr]{%
    1	4.99971\\
    2	4.99915\\
    3	5.00091\\
    4	5.00129\\
    5	5.00149\\
    6	4.99817\\
    7	5.0023\\
    8	4.99719\\
    9	5.00391\\
   10	4.99517\\
   11	4.99325\\
   12	5.00705\\
   13	4.99217\\
   14	5.00889\\
   15	4.98823\\
   16	4.98818\\
   17	4.98598\\
   18	5.01406\\
   19	4.98531\\
   20	4.98376\\
};
\addlegendentry{$\lambda_m$}

\addplot [color=red, draw=none, mark=o, mark options={solid, red}]
  table[row sep=crcr]{%
    1	4.99971\\
    2	4.99915\\
    3	5.00091\\
    5	5.00149\\
    6	4.99817\\
    7	5.0023\\
    8	4.99719\\
    9	5.00391\\
   10	4.99517\\
   11	4.99325\\
   12	5.00705\\
   13	4.99217\\
   14	5.00889\\
   15	4.98823\\
   16	4.98818\\
   17	4.98598\\
   19	4.98531\\
};
\addlegendentry{found $\lambda_m$}

\end{axis}
\end{tikzpicture}%
	\caption{}
\end{subfigure} \\
 \begin{subfigure}[b]{0.5\textwidth}
 	\hspace{0.7cm}
%
%
%
\definecolor{mycolor1}{rgb}{0.00000,0.44700,0.74100}%
\begin{tikzpicture}

\begin{axis}[%
width=0.951\figW,
height=0.75\figW,
at={(0\figW,0\figW)},
scale only axis,
xmin=    0,
xmax=   20,
ymin=  0.8,
ymax=  2.6,
ylabel style={font=\color{white!15!black}},
ylabel={$\lambda_m$},
axis background/.style={fill=white},
legend style={at={(0.03,0.97)}, anchor=north west, legend cell align=left, align=left, draw=white!15!black, only marks}
]
\addplot [color=mycolor1, draw=none, mark=asterisk, mark options={solid, mycolor1}]
  table[row sep=crcr]{%
    1	0.932368\\
    2	1.11385\\
    3	1.87369\\
    4	1.90999\\
    5	2.03706\\
    6	2.04216\\
    7	2.07979\\
    8	2.20542\\
    9	2.2147\\
   10	2.22965\\
   11	2.31916\\
   12	2.33797\\
   13	2.37888\\
   14	2.39569\\
   15	2.41483\\
   16	2.4229\\
   17	2.46996\\
   18	2.48392\\
   19	2.48858\\
   20	2.50796\\
};
\addlegendentry{$\lambda_m$}

\addplot [color=red, draw=none, mark=o, mark options={solid, red}]
  table[row sep=crcr]{%
    1	0.932368\\
    2	1.11385\\
    3	1.87369\\
    4	1.90999\\
    5	2.03706\\
    6	2.04216\\
    7	2.07979\\
    8	2.20542\\
    9	2.2147\\
   10	2.22965\\
   11	2.31916\\
   12	2.33797\\
   13	2.37888\\
   14	2.39569\\
   18	2.48392\\
   19	2.48858\\
};
\addlegendentry{found $\lambda_m$}

\end{axis}
\end{tikzpicture}%
	\caption{}
\end{subfigure}
\begin{subfigure}[b]{0.5\textwidth}
%
%
%
\definecolor{mycolor1}{rgb}{0.00000,0.44700,0.74100}%
\begin{tikzpicture}

\begin{axis}[%
width=0.951\figW,
height=0.75\figW,
at={(0\figW,0\figW)},
scale only axis,
xmin=    0,
xmax=   20,
ymin=4.985,
ymax= 5.01,
ylabel style={font=\color{white!15!black}},
ylabel={$\lambda_m$},
axis background/.style={fill=white},
legend style={at={(0.03,0.97)}, anchor=north west, legend cell align=left, align=left, draw=white!15!black, only marks}
]
\addplot [color=mycolor1, draw=none, mark=asterisk, mark options={solid, mycolor1}]
  table[row sep=crcr]{%
    1	5.0003\\
    2	4.99934\\
    3	4.99915\\
    4	4.9985\\
    5	4.99848\\
    6	4.9983\\
    7	5.00235\\
    8	5.00432\\
    9	4.9955\\
   10	4.99528\\
   11	4.99487\\
   12	5.00532\\
   13	5.00544\\
   14	4.99455\\
   15	4.99415\\
   16	5.00726\\
   17	5.00859\\
   18	4.99079\\
   19	5.00941\\
   20	4.98915\\
};
\addlegendentry{$\lambda_m$}

\addplot [color=red, draw=none, mark=o, mark options={solid, red}]
  table[row sep=crcr]{%
    4	4.9985\\
    6	4.9983\\
   11	4.99487\\
   13	5.00544\\
   14	4.99455\\
   15	4.99415\\
   16	5.00726\\
   18	4.99079\\
};
\addlegendentry{found $\lambda_m$}

\end{axis}
\end{tikzpicture}%
	\caption{}
\end{subfigure}
\caption{Figures (a)-(b) resp. (c)-(d) the results for the random $m$EP in Example \ref{fig:random_54EP}. We plot the wanted eigenvalue $\lambda_m$ in blue and the found eigenvalues in red. We see that if the wanted eigenvalues are exterior we can find them. In case they are interior or/and more close to each other, we see that it is far more difficult for the algorithm to find the wanted eigenvalues. Even after 100 sweeps we can only find a few.} \label{fig:random_54EP}
\end{figure}

\begin{eexample} \label{ex:random_54EP}
(Randomly generated 4EP and 5EP) We generate a random 4EP and 5EP with $n = 100$. We shift it such that all $\lambda_m$ are positive valued via \eqref{eqn:shift_eqns}. We perform two experiments, on the one hand we compute $\lambda_m$'s that are closest to $0$ and on the other hand we compute $\lambda_m$'s that are closest to $5$. In the former we perform 20 full sweeps and in the latter we perform 100 full sweeps. We see in Figure \ref{fig:random_54EP} (a) resp. (c) the $20$ smallest $|\lambda_m|$ of the 4EP resp. 5EP and we plot the $\lambda_m$'s that our algorithm found that were among these $20$ eigenvalues. In Figure \ref{fig:random_54EP} (b) resp. (d) the $20$ $\lambda_m$ closest to $5$ were plotted for the 4EP resp. 5EP together with $\lambda_m$'s that we found in our algorithm. We see that if the eigenvalues are exterior we can find almost all of the wanted values, if the eigenvalues are interior and more close to each other we find only a few.
\end{eexample}

\section{Dependency on parameters} \label{sect:NumExp_dependency}

In this section, we discuss how the computational time depends on $m$ and $n$ (size of the involved matrices). We divide the algorithm in four parts from which we do separate timings: explicit construction of $\VM X_{\neq k}^T \VM \Delta_0 \VM X_{\neq k}$, calculating eigenvalues of $\left( \VM X_{\neq k}^T \VM \Delta_m \VM X_{\neq k}, \VM X_{\neq k}^T \VM \Delta_m \VM X_{\neq k}\right)$, selecting and checking convergence (Algorithm~ \ref{algo:sel_eigv}) and update of the modes of $\VM X$ (see line \ref{line:reshape} and \ref{line:update} of Algorithm~\ref{algo:mEP_algo}).

\begin{figure}[h] 
\setlength{\figW}{6cm} 
 \begin{subfigure}[b]{0.5\textwidth}
 	\hspace{0.7cm}
%
%
%
\definecolor{mycolor1}{rgb}{0.00000,0.44700,0.74100}%
\definecolor{mycolor2}{rgb}{0.85000,0.32500,0.09800}%
\begin{tikzpicture}

\begin{axis}[%
width=0.951\figW,
height=0.75\figW,
at={(0\figW,0\figW)},
scale only axis,
xmin=    3,
xmax=   11,
xlabel style={font=\color{white!15!black}},
xlabel={$m$},
ymode=log,
ymin= 0.01,
ymax=13.4502,
yminorticks=true,
ylabel style={font=\color{white!15!black}},
ylabel={$t$ (in seconds)},
axis background/.style={fill=white},
title style={font=\bfseries},
title={Time for making $\VM X_{\neq k}^T \VM \Delta_0 \VM X_{\neq k}$},
legend style={at={(0.03,0.97)}, anchor=north west, legend cell align=left, align=left, draw=white!15!black}
]
\addplot [color=mycolor1]
  table[row sep=crcr]{%
    3	0.0398717\\
    4	0.0555105\\
    5	0.0790666\\
    6	0.127046\\
    7	0.217512\\
    8	0.51846\\
    9	1.39597\\
   10	3.781\\
   11	13.4502\\
};
\addlegendentry{Not reducing ranks}

\addplot [color=mycolor2]
  table[row sep=crcr]{%
    3	0.0348677\\
    4	0.0509752\\
    5	0.0784167\\
    6	0.130913\\
    7	0.203734\\
    8	0.347515\\
    9	0.587253\\
   10	0.761251\\
   11	1.1545\\
};
\addlegendentry{Reducing ranks}

\end{axis}
\end{tikzpicture}%
	\caption{}
\end{subfigure}
\begin{subfigure}[b]{0.5\textwidth}
	\hspace{0.7cm}  
%
%
%
\definecolor{mycolor1}{rgb}{0.00000,0.44700,0.74100}%
\definecolor{mycolor2}{rgb}{0.85000,0.32500,0.09800}%
\begin{tikzpicture}

\begin{axis}[%
width=0.951\figW,
height=0.75\figW,
at={(0\figW,0\figW)},
scale only axis,
xmin=    3,
xmax=   11,
xlabel style={font=\color{white!15!black}},
xlabel={$m$},
ymode=log,
ymin=    1,
ymax= 1000,
yminorticks=true,
ylabel style={font=\color{white!15!black}},
ylabel={$t$ (in seconds)},
axis background/.style={fill=white},
title style={font=\bfseries},
title={Time for finding $2b + q$ eigenvalues $\left( \VM X_{\neq k}^T \VM \Delta_m \VM X_{\neq k}, \VM X_{\neq k}^T \VM \Delta_m \VM X_{\neq k}\right)$},
legend style={at={(0.03,0.97)}, anchor=north west, legend cell align=left, align=left, draw=white!15!black}
]
\addplot [color=mycolor1]
  table[row sep=crcr]{%
    3	1.30893\\
    4	2.94035\\
    5	3.42421\\
    6	5.92292\\
    7	9.57393\\
    8	18.3998\\
    9	52.9966\\
   10	147.936\\
   11	555.664\\
};
\addlegendentry{Not reducing ranks}

\addplot [color=mycolor2]
  table[row sep=crcr]{%
    3	1.18762\\
    4	2.87551\\
    5	3.38164\\
    6	6.05147\\
    7	8.33044\\
    8	11.1107\\
    9	17.1417\\
   10	21.5296\\
   11	27.3496\\
};
\addlegendentry{Reducing ranks}

\end{axis}
\end{tikzpicture}%
	\caption{}
\end{subfigure} \\
 \begin{subfigure}[b]{0.5\textwidth}
 	\hspace{0.7cm}
%
%
%
\definecolor{mycolor1}{rgb}{0.00000,0.44700,0.74100}%
\definecolor{mycolor2}{rgb}{0.85000,0.32500,0.09800}%
\begin{tikzpicture}

\begin{axis}[%
width=0.951\figW,
height=0.75\figW,
at={(0\figW,0\figW)},
scale only axis,
xmin=    3,
xmax=   11,
xlabel style={font=\color{white!15!black}},
xlabel={$m$},
ymode=log,
ymin=   10,
ymax=100000,
yminorticks=true,
ylabel style={font=\color{white!15!black}},
ylabel={$t$ (in seconds)},
axis background/.style={fill=white},
title style={font=\bfseries},
title={Time for selecting $b$ eigenvalues and checking convergence},
legend style={at={(0.03,0.97)}, anchor=north west, legend cell align=left, align=left, draw=white!15!black}
]
\addplot [color=mycolor1]
  table[row sep=crcr]{%
    3	10.4851\\
    4	22.8014\\
    5	38.5239\\
    6	65.4146\\
    7	123.921\\
    8	315.405\\
    9	1165.15\\
   10	5424.44\\
   11	22416.6\\
};
\addlegendentry{Not reducing ranks}

\addplot [color=mycolor2]
  table[row sep=crcr]{%
    3	10.1865\\
    4	22.3834\\
    5	38.6738\\
    6	65.8309\\
    7	108.723\\
    8	180.974\\
    9	300.24\\
   10	464.503\\
   11	685.197\\
};
\addlegendentry{Reducing ranks}

\end{axis}
\end{tikzpicture}%
	\caption{}
\end{subfigure}
\begin{subfigure}[b]{0.5\textwidth}
	\hspace{0.7cm}  
%
%
%
\definecolor{mycolor1}{rgb}{0.00000,0.44700,0.74100}%
\definecolor{mycolor2}{rgb}{0.85000,0.32500,0.09800}%
\begin{tikzpicture}

\begin{axis}[%
width=0.951\figW,
height=0.75\figW,
at={(0\figW,0\figW)},
scale only axis,
xmin=    3,
xmax=   11,
xlabel style={font=\color{white!15!black}},
xlabel={$m$},
ymin=    0,
ymax= 0.12,
ylabel style={font=\color{white!15!black}},
ylabel={$t$ (in seconds)},
axis background/.style={fill=white},
title style={font=\bfseries},
title={Time for the update of the modes $\VM X$},
legend style={at={(0.03,0.97)}, anchor=north west, legend cell align=left, align=left, draw=white!15!black}
]
\addplot [color=mycolor1]
  table[row sep=crcr]{%
    3	0.0234342\\
    4	0.032984\\
    5	0.0381419\\
    6	0.0485587\\
    7	0.0580388\\
    8	0.0689285\\
    9	0.0802892\\
   10	0.091895\\
   11	0.115818\\
};
\addlegendentry{Not reducing ranks}

\addplot [color=mycolor2]
  table[row sep=crcr]{%
    3	0.0192119\\
    4	0.0269206\\
    5	0.0368773\\
    6	0.0492817\\
    7	0.0574384\\
    8	0.0688035\\
    9	0.0779162\\
   10	0.0896916\\
   11	0.107256\\
};
\addlegendentry{Reducing ranks}

\end{axis}
\end{tikzpicture}%
	\caption{}
\end{subfigure}
\caption{Figures corresponding to Example \ref{fig:increasing_m}. We plot the time needed to perform the different parts of the algorithm when we increase the number of parameters $m$. We see that the computational times decrease by almost a factor ten when reducing first the ranks of the TT-operator for $m = 11$. We see that solving the projected eigenvalue problem is most computationally demanding. Note that all but the last plot is semi-log scaled. } \label{fig:increasing_m}
\end{figure}
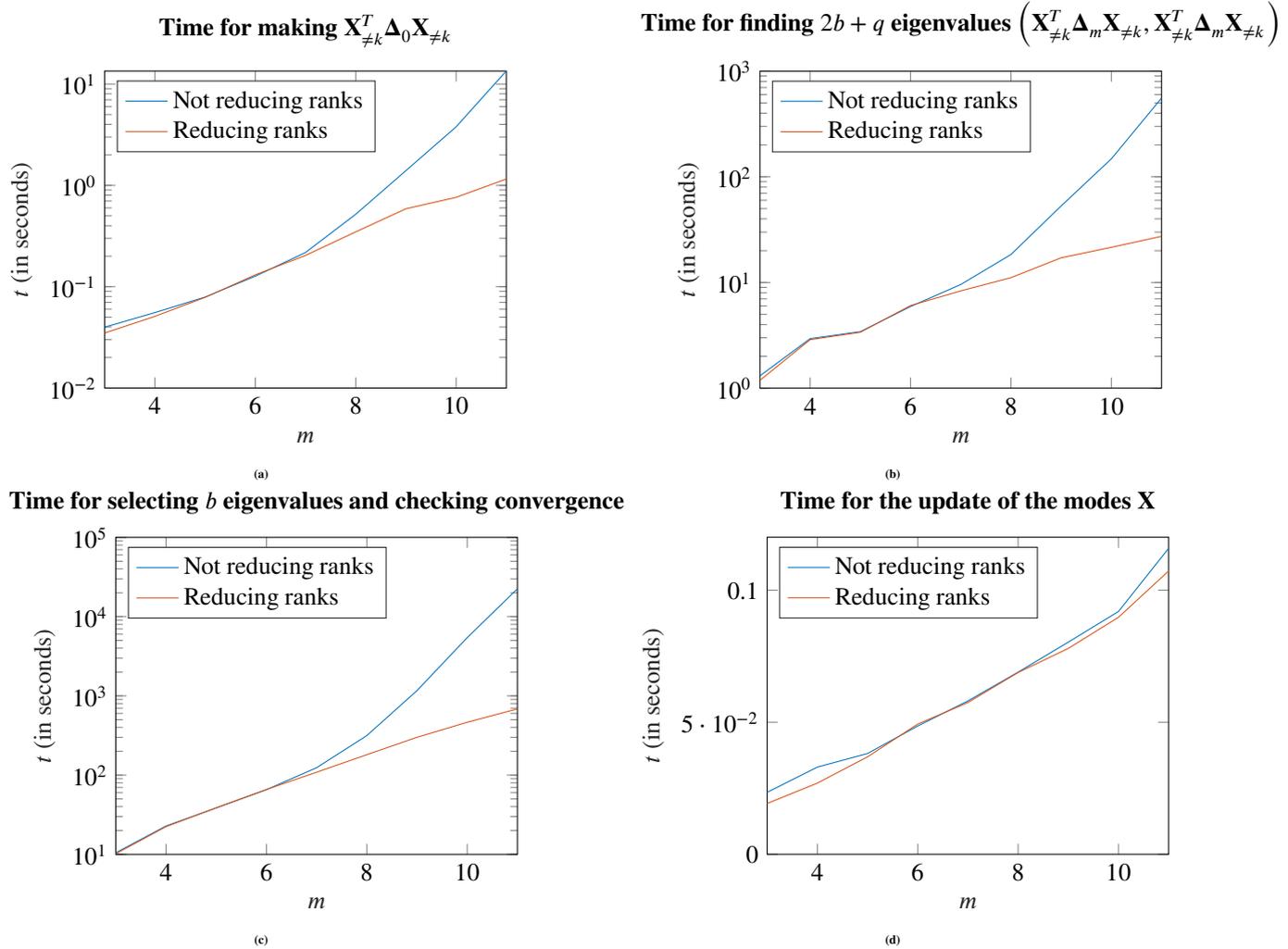

\begin{eexample} \label{ex:increasing_m}
(Increasing $m$) We fix $n=10$, make random $m$EP's for $m = 3,4,\hdots,11$ and we shift each $m$EP via \eqref{eqn:shift_eqns} such that the searched eigenvalues are exterior. We time the different parts as discussed above when performing $20$ full sweeps and we take the average over ten separate runs. The internal ranks of the $\VM \Delta$-matrices are the $m+1$st row of Pascal's triangle but we know that the actual ranks are maximal $n^2$. So for $m > 8$, the maximal internal rank can be reduced. To see what happens in the worst case (no reduction possible for the ranks), we do first an experiment without reducing the internal ranks of the $\VM \Delta$-matrices.
We perform the same experiment but now we reduce the internal rank by rounding the $\VM \Delta$-matrices as TT-operator to a tolerance of $10^{-13}.$  In Figure \ref{fig:increasing_m} we see the timings for the different parts when we round the $\VM \Delta$-matrices. 
\end{eexample}

We see in Figure~\ref{fig:increasing_m} (a)-(c) exponential dependency on $m$, when the $\VM \Delta$-matrices are not rounded. As well for the matrix-vector product $\VM \Delta_j \VM v$ or $\VM X_{\neq k}^T \VM \Delta_j \VM X_{\neq k} \VM v$ as well for the calculation of a projection $\VM X_{\neq k}^T \VM \Delta_0 \VM X_{\neq k}$, the number of operations depend on $m$ as $m (r^\Delta)^2$. On average $r^\Delta$ is $\dfrac{2^m}{m}$ as the ranks are the m+1st line in Pascal's triangle. This means that if $m$ is increased by one, we need to perform almost four times as much work. When we round the $\VM \Delta$-matrices first, we see that the computational time is reduced by a factor ten for $m = 11$. In the update of the modes of $\VM X$, the number of operations does only depend linearly on $m$ as no matrix-vector operations with one of the $\VM \Delta$-matrices is performed here. This explains Figure~\ref{fig:increasing_m}~(d) where we see that the rounding of the $\VM \Delta$-matrices does not affect the computational time for the update of the modes of $\VM X$.
 
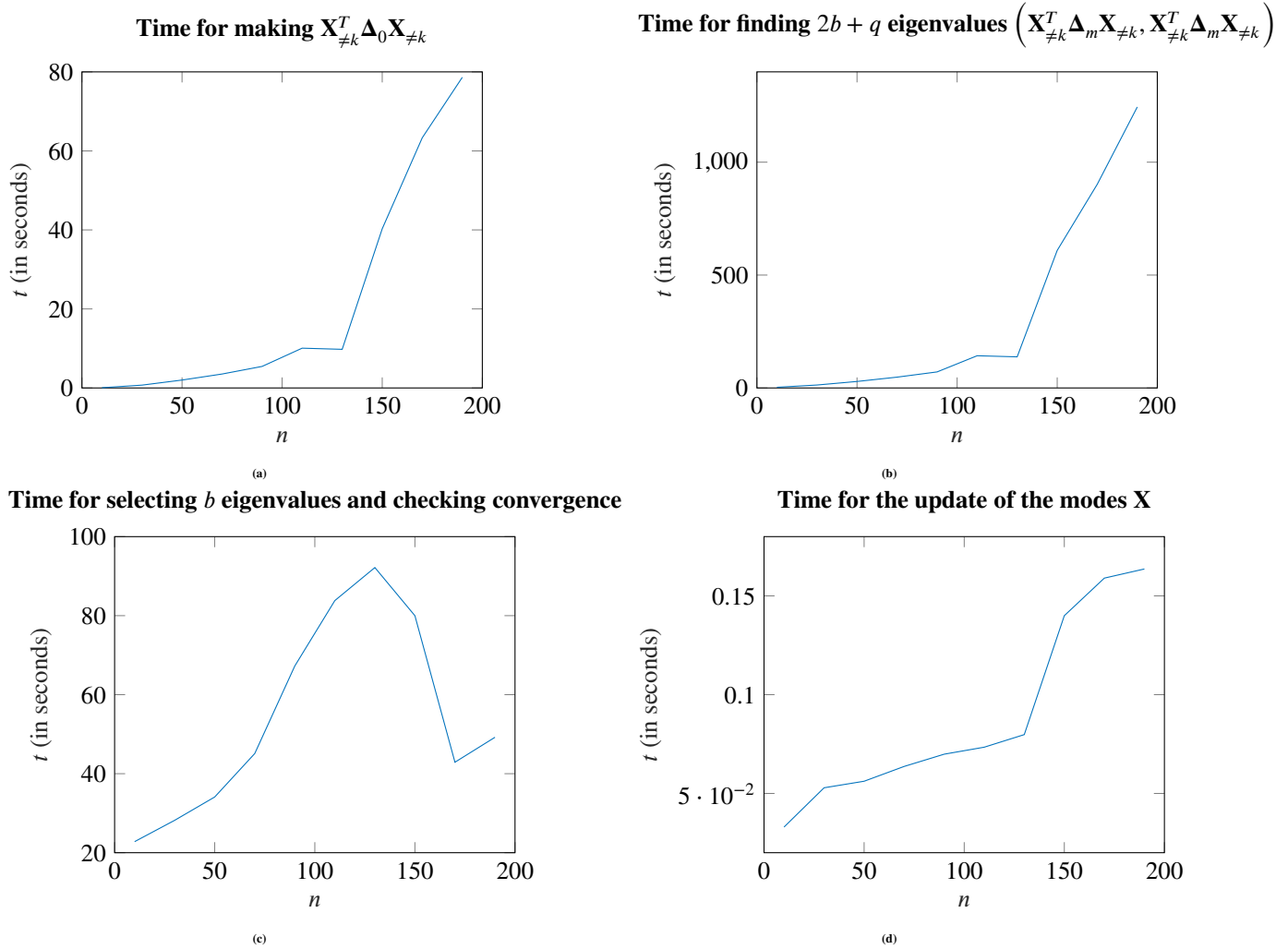
\begin{figure}[h] 
\setlength{\figW}{6cm} 
 \begin{subfigure}[b]{0.5\textwidth}
 	\hspace{0.7cm}
%
%
%
\definecolor{mycolor1}{rgb}{0.00000,0.44700,0.74100}%
\begin{tikzpicture}

\begin{axis}[%
width=0.951\figW,
height=0.75\figW,
at={(0\figW,0\figW)},
scale only axis,
xmin=    0,
xmax=  200,
xlabel style={font=\color{white!15!black}},
xlabel={$n$},
ymin=    0,
ymax=   80,
ylabel style={font=\color{white!15!black}},
ylabel={$t$ (in seconds)},
axis background/.style={fill=white},
title style={font=\bfseries},
title={Time for making $\VM X_{\neq k}^T \VM \Delta_0 \VM X_{\neq k}$}
]
\addplot [color=mycolor1, forget plot]
  table[row sep=crcr]{%
   10	0.0555105\\
   30	0.713559\\
   50	1.99786\\
   70	3.52555\\
   90	5.44928\\
  110	10.0675\\
  130	9.7764\\
  150	40.2867\\
  170	63.3273\\
  190	78.6256\\
};
\end{axis}
\end{tikzpicture}%
	\caption{}
\end{subfigure}
\begin{subfigure}[b]{0.5\textwidth}
	\hspace{0.7cm}  
%
%
%
\definecolor{mycolor1}{rgb}{0.00000,0.44700,0.74100}%
\begin{tikzpicture}

\begin{axis}[%
width=0.951\figW,
height=0.75\figW,
at={(0\figW,0\figW)},
scale only axis,
xmin=    0,
xmax=  200,
xlabel style={font=\color{white!15!black}},
xlabel={$n$},
ymin=    0,
ymax= 1400,
ylabel style={font=\color{white!15!black}},
ylabel={$t$ (in seconds)},
axis background/.style={fill=white},
title style={font=\bfseries},
title={Time for finding $2b + q$ eigenvalues $\left( \VM X_{\neq k}^T \VM \Delta_m \VM X_{\neq k}, \VM X_{\neq k}^T \VM \Delta_m \VM X_{\neq k}\right)$}
]
\addplot [color=mycolor1, forget plot]
  table[row sep=crcr]{%
   10	2.94035\\
   30	12.9417\\
   50	28.6961\\
   70	47.8382\\
   90	70.9725\\
  110	142.597\\
  130	137.755\\
  150	609.146\\
  170	902.038\\
  190	1244.44\\
};
\end{axis}
\end{tikzpicture}%
	\caption{}
\end{subfigure} \\
 \begin{subfigure}[b]{0.5\textwidth}
 	\hspace{0.7cm}
%
%
%
\definecolor{mycolor1}{rgb}{0.00000,0.44700,0.74100}%
\begin{tikzpicture}

\begin{axis}[%
width=0.951\figW,
height=0.75\figW,
at={(0\figW,0\figW)},
scale only axis,
xmin=    0,
xmax=  200,
xlabel style={font=\color{white!15!black}},
xlabel={$n$},
ymin=   20,
ymax=  100,
ylabel style={font=\color{white!15!black}},
ylabel={$t$ (in seconds)},
axis background/.style={fill=white},
title style={font=\bfseries},
title={Time for selecting $b$ eigenvalues and checking convergence}
]
\addplot [color=mycolor1, forget plot]
  table[row sep=crcr]{%
   10	22.8014\\
   30	28.1998\\
   50	34.1029\\
   70	45.1014\\
   90	67.2614\\
  110	83.8144\\
  130	92.1759\\
  150	79.9884\\
  170	42.9064\\
  190	49.2142\\
};
\end{axis}
\end{tikzpicture}%
	\caption{}
\end{subfigure}
\begin{subfigure}[b]{0.5\textwidth}
	\hspace{0.7cm} 
%
%
%
\definecolor{mycolor1}{rgb}{0.00000,0.44700,0.74100}%
\begin{tikzpicture}

\begin{axis}[%
width=0.951\figW,
height=0.75\figW,
at={(0\figW,0\figW)},
scale only axis,
xmin=    0,
xmax=  200,
xlabel style={font=\color{white!15!black}},
xlabel={$n$},
ymin= 0.02,
ymax= 0.18,
ylabel style={font=\color{white!15!black}},
ylabel={$t$ (in seconds)},
axis background/.style={fill=white},
title style={font=\bfseries},
title={Time for the update of the modes $\VM X$}
]
\addplot [color=mycolor1, forget plot]
  table[row sep=crcr]{%
   10	0.032984\\
   30	0.0528618\\
   50	0.0561787\\
   70	0.06367\\
   90	0.0698715\\
  110	0.0734119\\
  130	0.0797783\\
  150	0.139976\\
  170	0.158996\\
  190	0.163619\\
};
\end{axis}
\end{tikzpicture}%
	\caption{}
\end{subfigure}
\caption{Figures corresponding to Example \ref{fig:increasing_n}. We plot the time needed to perform the different parts of the algorithm when we increase the size $n$ of the matrices. We see that the solving of the projected eigenvalue problem is the most computationally demanding.  The computational time needed for the selection and the convergence, depend on how many eigenvalues are converged and at which sweep it converges, explaining the decrease in the computational time in (c). } \label{fig:increasing_n}
\end{figure}

\begin{eexample} \label{ex:increasing_n}
(Increasing $n$) We fix $m = 4$, make random $m$EP's for $n = 10, 30,\hdots, 190$ and we shift it such that the searched eigenvalues are exterior using \eqref{eqn:shift_eqns}. We time the different parts in the same way as in Example \ref{ex:increasing_m} when performing $20$ sweeps and we take the average over ten separate runs. The results are plotted in Figure~\ref{fig:increasing_n}. The computational time needed for the selection and the convergence, depend on how many eigenvalues are converged and at which sweep it converges, explaining the decrease in the computational time in Figure \ref{fig:increasing_n} (c). In all other parts we see a quadratic increase explained by the fact that there a parts for which the complexity depends quadratically on $n$. We see that the solving of the projected eigenvalue problems takes most of the time.
\end{eexample}

\section{Conclusions}
In this paper we presented a method that makes use of tensors to calculate a selection of the eigenvalues of a $m$EP with $m > 3$ and $n$ reasonably large. We showed that the Tensor-Train format can efficiently represent the $m$EP and showed that the associated $\VM \Delta$-matrices can be seen as Tensor-Train operators. This allows us to modify the algorithm in \cite{Dolgov2014a} to solve $m$EPs. Making modifications to the selection- and convergence criterion is necessary as we need to compute the eigenvalues of a non-symmetric generalised eigenvalue problem. We describe our modifications and motivate why we made them.

We illustrate with numerical experiments that our algorithm is able to find eigenvalues near zero and if they are well separated, we find all wanted eigenvalues.

\section*{Acknowledgements}
The authors thank the referees for their helpful remarks. This work was supported by the project KU Leuven Research Council grant C14/17/072 and by the project G0A5317N of the Research Foundation-Flanders (FWO - Vlaanderen). This work does not have any conflicts of interest.

\bibliography{references}%

\end{document}